\documentclass[12pt]{amsart}
\usepackage{subfigure}
\usepackage{graphicx}
\usepackage{rotating}
\usepackage{morefloats}
\usepackage[T1]{fontenc}
\usepackage[utf8]{inputenc}
\usepackage{amssymb}
\usepackage{amsmath}    
\def\R{I\kern -0,37 em R}
\def\P{I\kern -0,37 em P}
\def\Z{I\kern -0,37 em Z}

\title[Lie groupoids]{An introduction to Lie groupoids}
\author{A. Kumpera}
\address[Antonio Kumpera]{Campinas State University, Campinas, SP, Brazil}
\email{antoniokumpera@gmail.com}
\date{November, 2015}


\begin{document}

\maketitle

\begin{abstract}
We discuss the basic properties of Lie groupoids, Lie algebroids and Lie pseudo-groups in view of applying these techniques to the analysis of Jordan-Hölder resolutions and, subsequently, to the integration of partial differential equations. The present introduction is an extension to Lie groupoids, as far as possible, of the so well known properties and techniques much useful in Lie groups theory. We mention \textit{as far as possible} since, in the case of Lie groupoids, just the first Lie Theorem holds. As for the \textit{prolongation} algorithm, it is extremely useful when dealing with groupoids whereas rather senseless in the case of Lie groups. 
\end{abstract}

\section{Introduction}
Differential equations have been studied for a long time since Newton, Leibniz and Monge. Nevertheless, we owe first to Sophus Lie and later to Élie Cartan the structuring of this problem with the introduction of the \textit{finite and infinite continuous groups of transformations} leaving invariant these equations. It is however worthwhile to observe (\cite{Kumpera1980,Kumpera1999,Kumpera2015}) that whereas 
Lie was most incredibly successful in dealing with equations of finite type (the solutions depending only upon a finite number of parameters) for he had a deep knowledge and understanding of finite continuous groups (finite dimensional Lie groups), in the general context he was unable to go any further beyond examining a few specific examples involving the simple infinite and transitive continuous groups (transitive Lie pseudo-groups). For example, he showed that Jacobi’s last multiplier method was the best possible integration method - shortest, most accurate and of lowest degree - due to the fact that the infinite continuous group of all volume preserving transformations is simple. He was aware of the four classes of (complex) transitive simple pseudo-groups of transformations but rather uneasy whether these were the only ones, thus preventing him to take any benefits stemming from a systematic use of Jordan-Hölder resolutions (\cite{Kumpera1980}). Consequently, it remained to Cartan to develop the infinite dimensional theory with a touch only accessible to the most illuminated. In his doctoral dissertation  (\cite{Cartan1894}), Cartan classified the transitive complex simple Lie pseudo-groups showing that these were precisely the four classes announced by Sophus Lie (see also \cite{Cartan1909} and \cite{Kuranishi1962}). Much later (\cite{Cartan1913}), he obtained the classification of the real forms and could start definitely working on what really mattered to him namely, the devising of integration processes for differential systems (partial differential equations) with the help of Jordan-Hölder resolutions for the Lie pseudo-groups of all the local equivalences, defined on the base 
spaces, of the given systems (\cite{Cartan1904, Cartan1905, Cartan1908, Cartan1909, Cartan1921, Cartan1930, Cartan1938/39}).

\vspace{3 mm}
\noindent
Lie pseudo-groups are families of local transformations, given on a manifold \textit{M} and defined as the solutions of a system of partial differential equations $\mathcal{S}$ of a given order \textit{k}. It is further required that the usual group theoretic properties namely, the stability under composition and the passage to inverses, hold and the basic property states that such a family is a Lie pseudo-group if abd only if the corresponding differential equation $\mathcal{S}$ constitutes a Lie sub-groupoid of the groupoid $\Pi_kM$ of all the $k-$th order invertible jets. We perceive thereafter the fundamental importance of Lie groupoids in the study of Lie pseudo-groups and thereafter in the study of differential systems. It should also be mentioned that these structures were introduced by Ehresmann in the midst of the last century (\cite{Ehresmann1958}).

\vspace{3 mm}
\noindent
The differential properties of Lie groupoids are subordinate to an algebraic structure namely, that of a groupoid structure. However, the restrictions thus imposed are much weaker than in the case of Lie groups since, in general, there is a small part of the groupoid, namely the subset of units, which does not carry any algebraic structure at all, except a trivial action of the groupoid itself, and the composition in the groupoid is only partially defined. As in the case of Lie groups, many algebraic and differential properties of the Lie groupoids can be described as well as derived in terms of the corresponding algebraic properties of their Lie algebras - Lie algebroids in the present context - thus giving rise to functors from an algebraic to a differential category. The local results are analogous to those for Lie groups, their proof requiring just a little more patience and effort. Quite to the contrary, the global results are sensibly more sophisticated and certainly more difficult to obtain. We do not believe it to be necessary to mention, in more details, the various subjects that shadll be discuss hereafter since a little familiarity with Lie groups in association with the above keywords is entirely self explaining.

\vspace{3 mm}
\noindent
The importance of Lie groupoids relies more in their usefulness as a technique to study geometrical problems than as mathematical objects \textit{per se}. Moreover, they constitute at present the most appropriate setting for the development of certain geometric theories \textit{e.g.}, infinitesimal connexions and, of course, Lie pseudo-groups. They find inasmuch their place in the local and global equivalence problem for geometric structures (\cite{Kumpera2015}) and, most important, in the study of differential systems and partial differential equations. We finally observe that the study of G-structures can be considered as a special case of the study of transitive Lie groupoids and, in the bibliography, we mention several articles that apply Lie groupoids in various differential geometric situations.

\vspace{3 mm}
\noindent
Most of the definitions and concepts adopted here are standard and follow the common consensus. As for the terminology, we do not have any privileged choices and find that the best option results in looking at the not always so explicit origins. We finally remark that all the spaces and structures considered here are \textit{small} i.e., defined on underlying sets and never on \textit{classes} as considered, in a much broader setting, by Ch. Ehresmann (\cite{Ehresmann1958}). 

\section{Definitions}
A partially defined internal operation in a set $\mathcal{A}$ is a mapping $\mathcal{B}\longrightarrow\mathcal{A}$ where $\mathcal{B}\subset\mathcal{A}\times\mathcal{A}.$ As is usual, we indicate $\mathcal{B}$ by $\mathcal{A}*\mathcal{A},$ say that \textit{a} is composable with \textit{b} (\textit{a},\textit{b} in this order) if $(a,b)\in\mathcal{A}*\mathcal{A}$ and write \textit{ab} for their composite. An element $e\in\mathcal{A}$ is called a unit if $ae=a$ for any element $a\in\mathcal{A}$ that is composable with \textit{e} (in this order) and inasmuch $eb=b$. We could as well consider left and right units but these are completely irrelevant to us.

\newtheorem{struct}[DefinitionCounter]{Definition}
\begin{struct}
A groupoid is a non-empty set $\Gamma$ together with a partially defined internal operation $\Gamma*\Gamma\longrightarrow\Gamma$ satisfying the following properties:

\vspace{4 mm}
$1$. (associativity) If (g,h),(h,k)$\in\Gamma*\Gamma$ then (gh,k),(g,hk)$\in\Gamma*\Gamma$ and (gh)k=g(hk).

\vspace{3 mm}
$2$. (existence of units) For any g $\in\Gamma,$ there exists a unit e which is right composable with g (right unit) and a unit e' left composable with g (left unit) i.e., ge=e'g=g.

\vspace{3 mm}
$3$. (existence of inverses) For any g $\in\Gamma$, there exists an element $g^{-1}$ which is left and right composable with g and such that $gg^{-1}=e'$ and $g^{-1}g=e$ where e and e' are units.
\end{struct}

\noindent
From (1) and (3), we infer that $(gh,k)\in\Gamma*\Gamma$ if and only if $(h,k)\in\Gamma*\Gamma$. A simple computation shows that the left and the right units associated to an element \textit{g} as well as the inverse element $g^{-1}$ are  unique. We define the \textit{source} map $\alpha:\Gamma\longrightarrow\Gamma_0,~~\alpha(g)=e,$ as well as the \textit{target} map  $\beta:\Gamma\longrightarrow\Gamma_0,~~\beta(g)=e'.$ The maps $\alpha$ and $\beta$ are retracts ($\alpha^2=\alpha, \beta^2=\beta$) whose common image is the sub-set $\Gamma_0$ of all the units of $\Gamma.$ The element \textit{g} is composable with \textit{h} (in this order) if and only if $\alpha(g)=\beta(h)$ and the inversion map $\varphi:g\in\Gamma\longmapsto g^{-1}\in\Gamma$ is an involution ($\varphi^2=Id$). There are several equivalent ways of stating the groupoid axioms (\cite{Barros1969}), the associativity law being usually expressed by:
\begin{equation*}
(g,h),(gh,k)\in\Gamma*\Gamma\Longleftrightarrow (h,k),(g,hk)\in\Gamma*\Gamma
\end{equation*}
and, this being the case, the identity $(gh)k=g(hk)$ holds.

\newtheorem{ehres}[DefinitionCounter]{Definition}h
\begin{ehres}
A differentiable groupoid is a non-empty set $\Gamma$ together with a differentiable structure and a groupoid structure satisfying the following compatibility conditions:

\vspace{4 mm}
$1$. The set of units $\Gamma_0$ is a sub-manifold of $\Gamma.$

\vspace{3 mm}
$2$. The source and target maps $\alpha,\beta:\Gamma\longrightarrow\Gamma_0$ are differentiable and transverse.

\vspace{3 mm}
$3$. The map $(g,h)\in\Gamma\times_{\Gamma_0}\Gamma\longmapsto gh\in\Gamma$ is differentiable, where $\Gamma\times_{\Gamma_0}\Gamma=\Gamma*\Gamma$ is the fibre product of $\alpha$ and $\beta.$

\vspace{3 mm}
$4$. The map $g\in\Gamma\longmapsto g^{-1}\in\Gamma$ is differentiable.
\end{ehres}

\noindent
The natural inclusion $\iota:\Gamma_0\longrightarrow\Gamma$ is a right inverse to $\alpha$ i.e., $\alpha\circ\iota=Id_{\Gamma_0}$ and, inasmuch, is also a right inverse to $\beta$ hence $\alpha$ and $\beta$ are of maximal rank at every point of $\Gamma_0.$ Consequently, they are also of maximal rank in a whole neighborhood $\mathcal{U}$ of $\Gamma_0$ which we can choose to be symmetric. A simple argument shows further that $\alpha$ and $\beta$ are of maximal rank at every point of the sub-groupoid $\Gamma(\mathcal{U})$ of $\Gamma$ generated by $\mathcal{U}$ (all finite products and inverses of the element belonging to $\mathcal{U}$) which in fact will be an open sub-groupoid of $\Gamma.$ It follows that for any differentiable groupoid $\Gamma,$ the map $\alpha$ and, by inversion, the map $\beta$ are of maximal rank on the $\alpha-$connected component of $\Gamma_0$ (union of the connected components, in the $\alpha-$fibres, of the points in $\Gamma_0,$ assumed to be connected), as well as on the $\beta-$connected component. 

\newtheorem{komo}[DefinitionCounter]{Definition}
\begin{komo}
A Lie groupoid is a differentiable groupoid for which the map $\alpha$ is a fibration (surmersion). By  inversion, $\beta$ is equally a fibration. 
\end{komo}

\noindent
We infer that, for Lie groupoids, the condition (2) of the previous definition is a consequence of the present hypothesis made on $\alpha.$ Let $\Gamma$ be a Lie groupoid, denote by $\Gamma_{\alpha e}$ the fibre of $\alpha$ over the point $e\in\Gamma_0,$ called the $\alpha-$fibre, $\Gamma_{\beta e}$ the $\beta-$fibre, $\Gamma_e$ the \textit{isotropy group} at the unit point \textit{e} and let us proceed in defining the infinitesimal counterparts. For this, we consider any element $h\in\Gamma$ and define the diffeomorphism - right translation - $\phi_h:g\in\Gamma_{\alpha e'}\longmapsto gh\in\Gamma_{\alpha e},$ $e=\alpha(h),~e'=\beta(h).$ A local vector field $\xi$ defined on $\Gamma$ is said to be \textit{right-invariant} when it is $\alpha-vertical$ i.e., vanishes by $\alpha_*$, and  $(\phi_h)_*\xi=\xi$ for every $h\in\Gamma$ whenever this relation has a meaning. When $\xi$ is right-invariant and defined on an open set $\mathcal{U},$ it can be extended to a right-invariant vector field $\xi'$ defined on the open set $\mathcal{V}=\beta^{-1}(\beta(\mathcal{U})).$ In fact, if $h\in\mathcal{V}$ and $g\in\mathcal{U},$ with $\beta{h}=\beta{g},$ we set $\xi'_h=(\phi_{\ell})_*\xi_g$ where $\ell=g^{-1}h.$ By restriction of $\xi'$, the field $\xi$ determines an $\alpha-$vertical vector field $\xi_0$ along U$=\beta(\mathcal{U})\subset\Gamma_0$ and where $(\xi_0)_e=(\phi_{g^{-1}})_*\xi_g\in T_y\alpha^{-1}(y)$ for any $g\in\mathcal{U}$ and $e=\beta(g).$ Conversely, if $\xi_0$ is an $\alpha-$vertical vector field defined along the open set U in $\Gamma_0$, then there exists a unique right-invariant vector field $\xi$ defined on $\beta^{-1}(U)$ which induces $\xi_0.$ We infer that the natural domains of definition of the right-invariant vector fields on $\Gamma$ are the $\beta-$saturated open sets of $\Gamma.$ Let us denote by $V\Gamma$ the sub-bundle of $T\Gamma$ composed by the $\alpha-$vertical vectors (\textit{i.e.}, tangent to the $\alpha-$fibres). The previous discussion shows that the pre-sheaf of right invariant vector fields defined on $\beta-$saturated open sets of $\Gamma$ is canonically isomorphic to the pre-sheaf of local sections of the vector bundle $V\Gamma|_{\Gamma_0},$ restriction of $V\Gamma$ to the units sub-manifold $\Gamma_0.$ It is obvious that the bracket of two right-invariant vector fields is again right-invariant hence, by transport, we obtain a structure of $\bf{R}-$Lie algebras on the pre-sheaf $\Gamma_{loc}(V\Gamma|_{\Gamma_0}).$ Passing to germs, we obtain a sheaf of $\bf{R}-$Lie algebras $\mathcal{L}$ on the base space $\Gamma_0$ which is called the \textit{Lie algebroid} of $\Gamma.$ This bracket satisfies the identity
\begin{equation}
[f\mu,g\eta]=fg[\mu,\eta]+f[\vartheta(\beta_*\mu)g]\eta-g[\vartheta(\beta_*\eta)f]\mu
\end{equation}
\noindent
where $\mu,\eta\in\mathcal{L},$ and $f,g\in\mathcal{O}_{\Gamma_0},$ the sheaf of functions on $\Gamma_0,$ $\beta_*:\mathcal{L}\longrightarrow\underline{T\Gamma_0}$ is the extension, to germs, of $\beta_*:V\Gamma|_{\Gamma_0}\longrightarrow T\Gamma_0$ and $\vartheta(~)$ denotes the Lie derivative by vector fields.

\vspace{3 mm}
\noindent
Let us now say a few words on the isotropies of $\Gamma.$ Each isotropy, at a unit element, forms a group that has some differential properties though, in general, it does not possess any suitable Lie group structure unless $\Gamma$ is a Lie groupoid. It operates to the right on the fibre $\alpha^{-1}(e),$ the operation being induced by the product in $\Gamma.$ This right action is of principal type with respect to the map $\beta:\alpha^{-1}(e)\longrightarrow\Gamma_0$ i.e., the action is simply transitive on each $\beta-$fibre. For any $g\in\Gamma_e,$ the right translation $\tau_g$ on $\alpha^{-1}(e)$ is a diffeomorphism that preserves (leaves invariant) the $\beta-$fibres. Remark however that $\beta:\alpha^{-1}(e)\longrightarrow\Gamma_0$ is not necessarily of locally constant rank and that $\beta(\alpha^{-1}(e))$ needs not possess any suitable differentiable sub-manifold structure of $\Gamma_0$ (in such a way, for instance, that $\beta:\alpha^{-1}(e)\longrightarrow\beta(\alpha^{-1}(e))$ be a fibration). We can also repeat entirely the previous discussion with the projection $\beta$ and obtain inasmuch a second Lie algebroid that, however, is isomorphic to the previous one \textit{via} the inversion map $g\longrightarrow g^{-1}.$

\vspace{3 mm}
\noindent
Many of the elementary properties of Lie groups transcribe to Lie groupoids. For example, any local or global differentiable morphism of Lie groupoids $\varphi:\Gamma\longrightarrow\Sigma$ which, restricted to the units, is a diffeomorphism, induces an infinitesimal morphism $\varphi_*:\mathcal{L}\longrightarrow\mathcal{S}$ of Lie algebroids. Conversely, any such Lie algebroid morphism generates a local groupoid morphism defined in a neighborhood of the units sub-manifold. A proof of this statement, though in a more restricted context, can be found in \cite{Que1967} and \cite{Que1968}. Globalising such a local morphism $\varphi$ and associating to a given Lie algebroid $\mathcal{L},$ defined on a manifold \textit{M}, a Lie groupoid with an isomorphic algebroid are far more complex problems. They were studied initially by Pradines (\cite{Pradines1967},\cite{Pradines1968}) and much later placed in full display by Almeida and Molino (\cite{Almeida1985}), who exhibited a counter-example to Lie's third theorem, and still later by Mackenzie (\cite{Krill1987}). Fortunately, a positive result for Cartan pseudo-groups can be found in \cite{Kumpera1963} and later stated in full generality for (real analytic) Lie pseudo-groups by Kuranishi (\cite{Kuranishi1967}). Concerning Lie's second theorem, a counter-example was provided by the author in 1963 (\cite{Pradines1966}) and, up to now, there does not seem to be any effective criterion under which the statement of the theorem would hold. We shall confine our attention to the more elementary aspects of the theory since for our purposes we do not need any more than that.

\section{Sub-groupoids and quotients}
In this section we discuss very standard facts though some additional care must be taken in the case of groupoids. A Lie sub-groupoid of $\Gamma$ is of course a subset $\Sigma$ stable under the product operation and under the passage to inverses. Moreover it is assumed that $\Sigma$ is a (not necessarily regularly embedded) sub-manifold of $\Gamma$ and that the units space $\Sigma_0$ is a sub-manifold of $\Sigma.$ A \textit{morphism} $\varphi:\Gamma\longrightarrow\Gamma_1$ is a mapping that preserves the composition hence also preserves the inverses since forcibly it must preserve the units. Curiously enough, the image of a morphism is not always a Lie sub-groupoid of the image. On the one hand, the images of two non-composable elements can become composable and one cannot claim anything about their composite. Inasmuch, the image needs not be a sub-manifold. \textit{A priori}, the units space of a sub-groupoid needs not coincide with that of the ambient groupoid and we shall therefore assume right from the beginning the following

\vspace{3 mm}
\noindent
\textbf{Hypotheses:}

\vspace{3 mm}
1. The units space of a sub-groupoid coincides with the ambient units space.

\vspace{3 mm}
2. Any morphism becomes a diffeomorphism upon restriction to the units.

\vspace{3 mm}
\noindent
Under these conditions, the image of a morphism becomes of course a sub-groupoid. Let $\Gamma$ be a given Lie groupoid and let us take a Lie sub-groupoid $\Sigma$. Then clearly, every $\alpha-$fibre of $\Sigma$ is contained in the corresponding $\alpha-$fibre of $\Gamma.$ Next, for any given unit \textit{e}, recall that we denote by $\Gamma_{\alpha e}$ and $\Sigma_{\alpha e}$, the corresponding $\alpha-$fibres and by $\Gamma_e$, resp. $\Sigma_e,$ the corresponding isotropy groups at the unit point \textit{e}. The union of all the isotropy groups needs not be a differentiable sub-manifold of $\Gamma.$ However, these isotropies will become a locally trivial differentiable sub-bundle of, say, the $\alpha-$fibres or the $\beta-$fibres bundle, if and only if the fibre product mapping 
\begin{equation*}
\alpha\vee\beta:g\in\Gamma\longmapsto(\alpha(g),\beta(g))\in\Gamma_0\times\Gamma_0
\end{equation*}
is transverse to the diagonal in $\Gamma_0\times\Gamma_0.$ We shall say that $\Sigma$ is \textit{normal} or \textit{invariant} in $\Gamma$ whenever its isotropy bundle is \textit{normal} in $\Gamma.$ A little explanation seems to be required. What we are seeking for is the usual condition of normality namely, $\gamma\cdot\Sigma\cdot\gamma^{-1}\subset\Sigma$ for an arbitrary element $\gamma\in\Gamma.$ Since $\alpha(\gamma)=\beta(\gamma^{-1})$ and $\beta(\gamma)=\alpha(\gamma^{-1}),$ we infer that the previous condition is only applicable to isotropic elements of $\Sigma$ and that the resulting subset $\gamma\cdot\Sigma_e\cdot\gamma^{-1},~e=\alpha(\gamma),$ is the isotropy $\Sigma_{e'}$ of $\Sigma$ at the point $e'=\beta(\gamma).$ The above normality requirement means consequently that each subgroup $\Sigma_e$ is normal in $\Gamma_e$ and, furthermore, that the isotropy sub-bundle of $\Sigma$ is preserved by the (adjoint) action of arbitrary elements $\gamma\in\Gamma.$ By invertibility, the previous inclusion becomes in fact an equality.

\vspace{3 mm}
\noindent
We are now all set to define the quotient $\Gamma/\Sigma$ in the category of Lie groupoids. However, instead of defining the equivalence relation that will factor $\Gamma$ and ultimately produce the quotient, we shall first introduce the equivalence classes. Of course, we can consider either the \textit{right} classes $\Sigma\cdot\gamma=\Sigma_{\beta(\gamma)}\cdot\gamma$ or the \textit{left} classes $\gamma\cdot\Sigma=\gamma\cdot\Sigma_{\alpha(\gamma)}.$ On account of the inversion $\gamma\longmapsto\gamma^{-1}$, both procedures are equivalent and, on account of normality, both procedures will turn out to be equal. Since every $\alpha-$fibre of $\Sigma$ is contained in the corresponding $\alpha-$fibre of $\Gamma,$ we would like that each $\alpha-$fibre $\Sigma_{\alpha e}$ be an equivalence class that factors onto (identify with) the corresponding unit in the \textit{eventual} quotient groupoid. However, this cannot be so since the elements belonging to an $\alpha-$fibre can have many targets. The same problem appears with respect to the $\beta-$fibres. We shall therefore have to argue with the isotropy groups.

\vspace{3 mm}
\noindent
Let us, at present, just consider the right classes (right co-sets) and, for any element $\gamma\in\Gamma,$ we set $\Sigma_e\cdot\gamma,$ where $e=\beta(\gamma)$ ($\Sigma_e=\Sigma_e\cdot e$). We first show that two such classes either coincide or are disjoint. In fact, assuming that $Z$ belongs to both $\Sigma_e\cdot Y$ ($e=\beta(Y)$) and $\Sigma_{e'}\cdot Y'$, then $Z=X\cdot Y=X'\cdot Y'$, $\alpha(Z)=\alpha(Y)=\alpha(Y')$, $\Sigma_e\cdot Y\cdot(Y')^{-1}=\Sigma_{e'}$ hence ultimately $\Sigma_e \cdot Y=\Sigma_{e'}\cdot Y'$. Next, we define the \textit{equivalence relation} on $\Gamma$ whose equivalence classes are precisely the above co-sets. If $Z,~Z'$ both belong to the co-set $\Sigma_{\bf{e}}\cdot Y$, then $Z=X\cdot Y$ and $Z'=X'\cdot Y$ with $X,~X'\in\Sigma_{\bf{e}}$ whereupon $Z'\cdot Z^{-1}=X'\cdot X^{-1}\in\Sigma_e$. We infer that $Z\sim Z'$ if and only if $Z'\cdot Z^{-1}\in\Sigma,$ as in group theory. 

\vspace{3 mm}
\noindent
We next define the groupoid structure on the quotient $\Gamma/\Sigma,$ composed by the above co-sets, and assume that $\Sigma$ is topologically closed in $\Gamma$ otherwise the quotient topology will be inadequate to underlie a differentiable manifold structure. It follows thereafter that the Lie sub-group $\Sigma_e$ is also closed in $\Gamma_e.$ Setting $\alpha(\Sigma_e\cdot\gamma)=\alpha(\gamma),$ $\beta(\Sigma_e\cdot\gamma)=e$ and taking any two composable co-sets $\Sigma_e\cdot\gamma$ and $\Sigma_{\epsilon}\cdot\delta,$ $e=\beta(\gamma),~\epsilon=\beta(\delta),~\alpha(\gamma)=\beta(\delta),$ we define the product by 
\begin{equation}
(\Sigma_e\cdot\gamma)\cdot(\Sigma_{\epsilon}\cdot\delta)=\Sigma_e\cdot\gamma\cdot\delta.
\end{equation}
We now proceed to show that this operation does not depend upon the representatives. Assuming that $\Sigma_e\cdot\gamma=\Sigma_{e'}\cdot\gamma',$ then $\gamma'=X\cdot\gamma,$ with $X\in\Sigma_e,$ and consequently $\Sigma_{e'}\cdot\gamma'=\Sigma_e\cdot\gamma'$ hence,
\begin{equation*}
(\Sigma_e\cdot\gamma')\cdot(\Sigma_{\epsilon}\cdot\delta)=\Sigma_e\cdot\gamma'\cdot\delta=\Sigma_e\cdot X\cdot\gamma\cdot\delta=\Sigma_e\cdot\gamma\cdot\delta
\end{equation*}
since $\Sigma_e\cdot X=\Sigma_e.$ It should be observed that this first step does not require the normality condition on $\Sigma.$ Let us next assume that $\Sigma_{\epsilon}\cdot\delta=\Sigma_{\epsilon'}\cdot\delta'.$ Then $\delta'=X\cdot\delta,~X\in\Sigma_{\epsilon}$ and consequently
\begin{equation*}
(\Sigma_e\cdot\gamma)\cdot(\Sigma_{\epsilon'}\cdot\delta')=\Sigma_e\cdot\gamma\cdot\delta'=\Sigma_e\cdot\gamma\cdot X\cdot\delta=\Sigma_e\cdot\gamma\cdot X\cdot\gamma^{-1}\cdot\gamma\cdot\delta=\Sigma_e\cdot Y\cdot\gamma\cdot\delta,
\end{equation*}
where $Y=\gamma\cdot X\cdot\gamma^{-1}\in\Sigma_e,$ hence the above product is equal to $\Sigma_e\cdot\gamma\cdot\delta$ since $\Sigma_e\cdot Y=\Sigma_e.$ Needless to say that the inverse of $\Sigma_e\cdot\gamma$ is equal to $\Sigma_{\epsilon}\cdot\gamma^{-1}$ where $\epsilon=\alpha(\gamma).$ 

\vspace{3 mm}
\noindent
The fact that $\beta\vee\alpha$ is transverse to the diagonal in $M\times M$ guarantees that the quotient $\Gamma/\Sigma$ has a differentiable manifold structure compatible with its groupoid structure. We shall however skip the details and only remark that the above definition of the quotient does not agree with that given in \cite{Krill1987} and outlined in \cite{Pradines1966} though both are closely related since, apparently, one cannot escape from operating with isotropies. The advantage of our approach relies in the simplicity for expressing the equivalence class of a given element $\gamma\in\Gamma,$ identical to that found in group theory.

\vspace{3 mm}
\noindent
As for the infinitesimal aspects of the above discussion, there is not much more to be added. It suffices to define the isotropy sub-algebra of $\mathcal{L}$ as being the sub-sheaf composed by those germs whose representatives vanish at the unit element involved. We shall say that $\mathcal{S}$ is a sub-algebroid 
in ideals of $\mathcal{L}$ when each stalk of the isotropy of $\mathcal{S}$ is an ideal in the corresponding stalk of $\mathcal{L}$ at that point. Integrating this condition for each germ of vector field belonging to $\mathcal{L},$ we infer that $\mathcal{S}$ is invariant under the local actions associated to arbitrary sections of $\mathcal{L}.$ The process of defining now a quotient Lie algebroid becomes obvious. It is also self-evident that each isotropy algebra is the Lie algebra of the corresponding isotropy group and that the quotient algebroid is the Lie algebroid associated to the quotient Lie groupoid whenever such an association exists.

\section{The Exponential map}
Let $\Gamma$ be a Lie groupoid and $\mathcal{L}$ its algebroid. Since the latter is the sheaf of germs of local sections of $V\Gamma|_{\Gamma_0},$ then $\Gamma(U,\mathcal{L})=\Gamma(U,V\Gamma|_{\Gamma_0}),$ where \textit{U} is an open set in $\Gamma_0$ and where, at present, the notation $\Gamma(~)$ also indicates a pre-sheaf of local sections. A local section $\sigma\in\Gamma(U,\Gamma),$ with respect to $\alpha,$ is called \textit{admissible} if $V=\beta\circ\sigma(U)$ is open in $\Gamma_0$ and if moreover the mapping $\beta\circ\sigma:U\longrightarrow V$ is a diffeomorphism. The groupoid structure of $\Gamma$ extends to the pre-sheaf $\Gamma_{a,loc}(\Gamma_0,\Gamma)$ of all the admissible local sections. In fact, if $\sigma$ is defined on the open set \textit{U} and $\tau$ on the open set \textit{V}, with $\beta\circ\sigma(U)\subset V,$ we define $\tau\cdot\sigma$ by $\tau\cdot\sigma(e)=[\tau\circ\beta\circ\sigma(e)]\cdot\sigma(e)$ where the right hand side product is the composition in $\Gamma.$

\vspace{3 mm}
\noindent
We next consider a local section $\Xi\in\Gamma(\Gamma_0,\mathcal{L})$ of $\mathcal{L},$ defined on an open set \textit{U}, let $\xi$ be the corresponding right-invariant vector field on $\beta^{-1}(U)\subset\Gamma$ and $(\varphi_t),~\varphi_t=exp\!~t\xi,$ the associated local 1-parameter group of transformations generated by $\xi.$ If we set $Exp\!~t\Xi=\varphi_t\circ\iota:U\longrightarrow\Gamma$ with $\iota:\Gamma_0\longrightarrow\Gamma$ the inclusion then, for any $e\in U\subset\Gamma_0,$ there is an open neighborhood $U_e$ of \textit{e} and a constant  $\epsilon_e>0$ such that $Exp\!~t\Xi:U_e\longrightarrow\Gamma$ is defined on $U_e$ for $t<\epsilon_e.$ Since $\xi$ is tangent to the $\alpha-$fibres of $\Gamma,$ its integral curves are always contained in these fibres hence $(\varphi_t)$ also preserves them. It follows that $Exp\!~t\Xi:U_e\longrightarrow\Gamma$ is an $\alpha-$section for each fixed value of \textit{t}. Moreover, since $\xi$ is $\beta-$projectable, we infer that
\begin{equation*}
\beta\circ Exp\!~t\Xi|U_e=exp\!~t(\beta_*\xi)|U_e,
\end{equation*}
hence $Exp\!~t\xi|U_e$ is, for each fixed \textit{t}, an admissible local section of $\Gamma.$ Inasmuch, $Exp\!~t\xi|U_e$ is a differentiable local 1-parameter sub-group of the groupoid $\Gamma_{a,loc}(\Gamma_0,\Gamma)$ generated by $\Xi$ in the sense that

\vspace{3 mm}
1. Each local family $(Exp\!~t\Xi|_{U_e}),$ $|t|<\epsilon_e,$ depends differentiably on the parameter \textit{t} and two such local families agree, for common values of \textit{t}, on the overlap of their domains,

\vspace{3 mm}
2. $\frac{d}{dt} Exp\!~t\Xi|_{t=0}=\Xi,$

\vspace{3 mm}
3. (a) $Exp\!~0=\iota,$ the units inclusion, and
 
\vspace{2 mm}
\hspace{4 mm}
(b) $Exp\!~(t+u)\Xi(e)=Exp\!~t\Xi(e')\cdot Exp\!~u\Xi(e),~e'=\beta\circ Exp\!~u\Xi(e),$ whenever both members are defined. 

\vspace{3 mm}
\noindent
We now wish to find conditions on $\Xi$ under which the Exponential map is global \textit{i.e.}, it is defined for every $t\in\bf{R}$ and is an admissible section defined on $\Gamma_0.$ Before stating the first Lemma, we recall that any right-invariant vector field, operating on a Lie groupoid $\Gamma,$ is always $\beta-$projectable and we shall  say that a given vector field, defined on a manifold \textit{M}, is \textit{global} when it generates a global 1-parameter group \textit{i.e.}, determined for every value of \textit{t}. Of course, the vector field is global if and only if its integral curves can be extended to infinity both ways.

\newtheorem{shit}[LemmaCounter]{Lemma}
\begin{shit}
Let $\Gamma$ be a Lie groupoid, $\xi$ a right-invariant vector field and $\zeta=\beta_*\xi$ the projected field on the units space $\Gamma_0.$ Then $\xi$ is 
global if and only if $\zeta$ is global. 
\end{shit}

\noindent
\textbf{Proof.} When $\xi$ is global then so is $\zeta$ since the integral curves of the latter are the $\beta-$projections of the integral curves of the former. For each fixed \textit{t}, the transformation $exp\!~t\xi$ is $\beta-$projectable and its projection, in $\Gamma_0,$ is equal to $exp\!~t\zeta.$ Assume conversely that $\zeta$ is global. Since $\xi$ is $\alpha-$vertical, any integral curve of $\xi$ is contained in an $\alpha-$fibre. Moreover, since $\xi$ is right-invariant, its integral curves are permuted by the right translations of $\Gamma.$ We infer that for any $h\in\Gamma$ with $\beta(h)=\epsilon\in\Gamma_0$ and any integral curve $\gamma(t)$ of $\xi$ with initial data $\gamma(0)=\epsilon,$ the right translated curve $\rho(t)=\gamma(t)\cdot h$ is an integral curve of $\xi$ with initial data $\rho(0)=h.$ To prove that $\xi$ is global, it suffices to show that for any unit $e\in\Gamma_0$ there is a global integral curve $\gamma:\textbf{R}\longrightarrow\alpha^{-1}(e)$ of $\xi$ with initial data $\gamma(0)=e.$ We distinguish two cases:

\vspace{3 mm}
(a)\hspace{4 mm} $\zeta_e=0.$

\vspace{2 mm}
\noindent
We observe that any integral curve $\gamma(t)$ of $\xi$ with initial value $\gamma(0)\in\Gamma_e$ (the isotropy group at \textit{e}) will be entirely contained in $\Gamma_e$ since $\rho(t)=\beta\circ\gamma(t)$ is an integral curve of $\zeta$ with initial data $\rho(0)=e$ and, since $\zeta_e=0,$ we infer that $\rho(t)=e$ for all the value of \textit{t}. Let $\gamma:~]\!\!-\!\epsilon,\epsilon[\longrightarrow\Gamma_e$ be an integral curve of $\xi$ with $\gamma(0)=e,$ the unit element of the isotropy group that coincides with the unit element of the groupoid, and take any element $h\in\Gamma_e.$ Then $\tilde{\gamma}(t)=\gamma(t)\cdot h$ is an integral curve of $\xi$ with the initial condition $\tilde{\gamma}(0)=h$ and defined on the same interval $]\!-\epsilon,\epsilon[.$ A standard continuation argument will then show that $\gamma$ can be extended to a global integral curve $\textbf{R}\longrightarrow\Gamma_e$ of $\xi.$ Perhaps it would have been more clever to observe, instead, that in this case $\xi_e\in T\Gamma_e$ and consequently that $\xi|_{\Gamma_e}$ is a right invariant vector field defined on this Lie group, hence global.  

\vspace{3 mm}
(b)\hspace{4 mm} $\zeta_e\neq 0.$

\vspace{2 mm}
\noindent 
Let $\rho:\textbf{R}\longrightarrow\Gamma_0$ be the global integral curve of $\zeta$ with initial data $\rho(0)=e.$ Then $\zeta_{\rho(t)}=\frac{d\rho}{dt}\neq 0$ for any $t\in\bf{R},$ hence $\rho$ is an immersion and its image is a one-dimensional sub-manifold of $\Gamma_0.$ Let $\gamma:~]a,b[\longrightarrow\alpha^{-1}(e)$ be the maximal integral curve of $\xi$ with initial value $\gamma(0)=e.$ We want to show that $]a,b[=\bf{R}.$ Assume, for example, that \textit{b} is finite. Then $\beta\circ\gamma$ is equal to the integral curve $\rho:~]a,b[\longrightarrow\Gamma_0$ and, taking a point $h\in\alpha^{-1}(e)$ for which $\beta(h)=\rho(b)$ as well as letting
\begin{equation*}
\delta_h:~]b-\epsilon,b+\epsilon[\longrightarrow\alpha^{-1}(e)
\end{equation*}
be an integral curve of $\xi$ with initial data $\delta_h(b)=h,$ we arrive at the following situation. The curve  $\beta\circ\delta_h$ is equal to $\rho:~]b-\epsilon,b+\epsilon[\longrightarrow\Gamma_0$ and, since $\zeta_{\rho(b)}\neq 0,$ this curve is injective - an integral curve does not cross with itself -. Moreover, since $\xi$ is right-invariant, the right action of $\Gamma_e$ on $\alpha^{-1}(e)$ permutes the integral curves of $\xi$ hence, for any $g\in\Gamma_e,$ the curve
\begin{equation*}
\delta_{hg}:~]b-\epsilon,b+\epsilon[\longrightarrow\alpha^{-1}(e),\hspace{5 mm}\delta_{hg}(t)=\delta_h(t)\cdot g,
\end{equation*}
is the integral curve of $\xi$ with initial data $\delta_{hg}(b)=hg$ and with the $\beta-$projection
\begin{equation*}
\beta\circ\delta_{hg}=\rho:~]b-\epsilon,b+\epsilon[\longrightarrow\Gamma_0.
\end{equation*}
All the curves $\delta_{hg}$ are injective since their projections are injective, the image sub-space $M=\rho(]b-\epsilon,b+\epsilon[)$ is a one-dimensional sub-manifold of $\Gamma_0$ and each sub-space $M_g=\delta_{hg}(]b-\epsilon,b+\epsilon[)$ is also a one-dimensional sub-manifold of $\alpha^{-1}(e).$ If $g,g'\in\Gamma_e$ and $g\neq g',$ then $M_g\cap M_{g'}=\varnothing.$ In fact, assuming that $w=\delta_{hg}(t)=\delta_{hg'}(u),$ then $\beta(w)=\rho(t)=\rho(u)$ and consequently, by injectivity, $t=u$ and $\delta_{hg}(t)=\delta_{hg'}(t).$ The local uniqueness of the integral curves implies that $\delta_{hg}=\delta_{hg'}$ hence  $hg=\delta_{hg}(b)=\delta_{hg'}(b)=hg'$ i.e., $g=h.$ We infer that the subset $A=\{w\in\alpha^{-1}(e)\!~|\!~\beta(w)\in M\}$ is the union of the disjoint sub-manifolds $M_g,~g\in\Gamma_e.$ Taking a value \textit{c} with $b-\epsilon<c<b,$ then $w=\gamma(c)\in A$ and consequently there exists a unique element $g\in\Gamma_e$ such that $\gamma(c)=\delta_{hg}(d)$ with $b-\epsilon<d<b+\epsilon.$ On the other hand, $\rho(c)=\beta\circ\gamma(c)=\beta\circ\delta_{hg}(d)=\rho(d)$ hence $c=d$ and therefore the two integral curves $\gamma$ and $\delta_{hg}$ have the same value at the point \textit{c}. We thus conclude that $\gamma\!~|\!~]b-\epsilon,b[=\delta_{hg}\!~|\!~]b-\epsilon,b[$ hence $\delta_{hg}$ is a continuation of $\gamma$ which thus ceases to be maximal.
        
\vspace{3 mm}
\noindent
\textbf{Remark.} The method of proof extends to the following more general situation. Let $\xi$ and $\zeta$ be as in the Lemma. Then, for any element $g\in\Gamma,$ the maximal integral curve $\gamma$ of $\xi$ with initial data $\gamma(0)=g$ projects, via $\beta,$ onto the maximal integral curve $\rho$ of $\zeta$ with the initial data $\rho(0)=\beta(g).$

\vspace{3 mm}
\noindent
The map $\beta:\Gamma\longrightarrow\Gamma_0$ induces, by differentiation, a linear map $\beta_*:V\Gamma|_{\Gamma_0}\longrightarrow T\Gamma_0$ hence an $\mathcal{O}_{\Gamma_0}-$linear sheaf map $\beta_*:\mathcal{L}\longrightarrow \underline{T\Gamma_0}.$

\newtheorem{expo}[CorollaryCounter]{Corollary}
\begin{expo}        
Let $\Gamma$ be a Lie groupoid with Lie algebroid $\mathcal{L}$ and let $\Xi\in\Gamma(\Gamma_0,\mathcal{L}).$ Then $\Xi$ is global (i.e., $Exp\!~t\Xi$ is global) if and only if the same is true for the vector field $\beta_*\circ\Xi$ (compare with \cite{Kumpera1972}, the Lemma of section 19, where uniform means global). 
\end{expo}

\noindent
\textbf{Proof.} When $\Xi$ is global then so is $\beta_*\circ\Xi$ for $exp\!~t(\beta_*\circ\Xi)=\beta\circ Exp\!~t\Xi.$ Conversely, if $\beta_*\circ\Xi$ is global then the right-invariant vector field $\xi$ associated to $\Xi$ is global hence $Exp\!~t\Xi= (exp\!~t\xi)\circ\iota$ is global.

\vspace{3 mm}
\noindent
Let $\Gamma_c(\Gamma_0,\mathcal{L})$ be the set of all local sections $\Xi$ of $\mathcal{L}$ such that $\beta_*\circ\Xi$ has compact support. Then the Lie algebroid structure of $\mathcal{L}$ extends to $\Gamma_c(\Gamma_0,\mathcal{L})$.

\newtheorem{exten}[CorollaryCounter]{Corollary}
\begin{exten}
Any $\Xi\in\Gamma_c(\Gamma_0,\mathcal{L})$ is global. 
\end{exten}

\noindent
\textbf{Proof.} The vector field $\theta=\beta_*\circ\Xi$ has compact support hence is global. We remark however that the right-invariant vector field $\xi$ associated to $\Xi$ needs not be compactly supported.

\vspace{3 mm}
\noindent
Let $\Gamma_a(\Gamma_0,\Gamma)$ be the set of admissible local sections of $\Gamma$ \textit{i.e.}, sections $\sigma$ with respect to the fibration $\alpha$ such that $\beta\circ\sigma$ is a local diffeomorphism of $\Gamma_0.$ The groupoid structure of $\Gamma$ (or of $\Gamma_{a,loc}(\Gamma_0,\Gamma)$) extends to a groupoid structure on $\Gamma_a(\Gamma_0,\Gamma)$ and enables us to define the Exponential map
\begin{equation}
Exp:\Gamma_c(\Gamma_0,\mathcal{L})\longrightarrow\Gamma_a(\Gamma_0,\Gamma)
\end{equation}
by $Exp\!~\Xi=(exp\!~\xi)\circ\iota,$ where $\xi$ is the right-invariant vector field associated to $\Xi$ and $exp\!~\xi=exp\!~1\xi.$ It is clear that $Exp\!~(t\Xi)=Exp\!~t\Xi$ since $exp\!~t\xi=exp\!~1\eta=exp(t\!~\xi)$ with $\eta=t\xi.$ This map is differentiable and satisfies the following properties:

\vspace{4 mm}
1. For each $\Xi\in\Gamma_c(\Gamma_0,\mathcal{L}),$ the map $t\in\textbf{R}\longmapsto Exp\!~t\Xi\in\Gamma_a(\Gamma_0,\Gamma)$ is a differentiable 1-parameter subgroup of $\Gamma_a(\Gamma_0,\Gamma)$ generated by $\Xi$ \textit{i.e.}, the following identities hold: $Exp\!~(t+u)\Xi=(Exp\!~t\Xi)\cdot(Exp\!~u\Xi),$ $Exp\!~0=\iota,$ $Exp-\!t\Xi=(Exp\!~t\Xi)^{-1}$ and finally $\frac{d}{dt}Exp\!~t\Xi|_{t=0}=\Xi,$

\vspace{3 mm}
2. $\beta\circ Exp\!~t\Xi=exp\!~t(\beta_*\circ\Xi)$ and

\vspace{3 mm}
3. $Exp$ is uniquely determined by the property (1).

\vspace{3 mm}
\noindent
As for the differentiability, it means that $Exp$ transforms differentiable families of sections of $\mathcal{L}$ into differentiable families of sections of $\Gamma.$ The uniqueness property is a consequence of the local uniqueness for the solutions of ordinary differential equations. Taking appropriate local coordinates, we can construct, locally, a system of ordinary differential equations whose solutions, with certain initial data, are (pieces of) the curves $Exp\!~t\Xi(e).$ The above mentioned properties are also trivial consequences of the corresponding properties for $exp\!~t\xi,$ where $\xi$ is a right-invariant vector field on $\Gamma.$

\vspace{3 mm}
\noindent
When $\Gamma$ is a Lie group (a Lie groupoid with a single unit), it is clear that \textit{Exp} is the usual exponential map \textit{exp}.

\vspace{3 mm}
\noindent
We end this section by examining some further analogies between Lie groups and Lie groupoids.

\vspace{3 mm}
\noindent
Let $\Gamma$ be a Lie groupoid. A diffeomorphism $\phi:\Gamma\longrightarrow\Gamma$ is called a right translation of $\Gamma$ if

\vspace{3 mm}
(a) $\phi$ preserves the fibration $\alpha:\Gamma\longrightarrow\Gamma_0$ and therefore induces a diffeomorphism $f:\Gamma_0\longrightarrow\Gamma_0$ on the base space of units such that the diagram

\vspace{2 mm}
\begin{equation*}
\Gamma~\xrightarrow{~\phi~}~\Gamma
\end{equation*}
\begin{equation*}
\alpha\downarrow\hspace{12 mm}\downarrow\alpha
\end{equation*}
\begin{equation*}
\Gamma_0\xrightarrow{~f~}\Gamma_0
\end{equation*}

\vspace{3 mm}
\noindent
commutes and

\vspace{2 mm}
(b) For any $e\in\Gamma_0,$ the restriction $\phi:\Gamma_{\alpha e}\longrightarrow\Gamma_{\alpha\epsilon},~\epsilon=f(e),$ is a right translation \textit{i.e.}, there exists an element $h(e)\in\Gamma_{\alpha\epsilon} (=\alpha^{-1}(\epsilon))$ such that $\phi(g)=\phi_{h(e)}(g)=g\cdot h(e),$ where $g\in\Gamma_{\alpha e}.$ This implies, in particular, that $\phi$ commutes with $\beta.$

\vspace{3 mm}
\noindent
Given a right translation $\phi:\Gamma\longrightarrow\Gamma$ then, for any $e\in\Gamma_0,$ the element $h(e)$ of the condition (b) is unique, $h(e)=\phi(e)$ and $\beta[h(e)]=e.$ The mapping $h:e\in\Gamma_0\longmapsto h(e)\in\Gamma$ is an $\alpha-$admissible section of $\beta$ since $\alpha\circ h=f$ is a diffeomorphism of $\Gamma_0.$ Equivalently, the map $\epsilon\in\Gamma_0\longmapsto h\circ f^{-1}(\epsilon)$ is a $\beta-$admissible section of $\alpha$ which, by inversion, yields the $\beta-$admissible section $\sigma:e\in\Gamma_0\longmapsto h(e)^{-1}\in\Gamma$ of the projection $\alpha,$ where $h(e)^{-1}$ is the inverse of $h(e)$ in $\Gamma.$ The right translation $\phi$ is given by
\begin{equation}
\phi(g)=g\cdot\sigma(e)^{-1},\hspace{5 mm} e=\alpha(g).
\end{equation}

\vspace{3 mm}
\noindent
Conversely, if $\sigma\in\Gamma_a(\Gamma_0,\Gamma),$ then the relation (4) defines a right translation that induces $\sigma.$ We infer that the group $\Gamma_a(\Gamma_0,\Gamma)$ of global sections identifies canonically with the group of right translations of $\Gamma$ and this identification is a group isomorphism.

\vspace{5 mm}
\noindent
Interchanging $\alpha$ with $\beta,$ we can define left translations of $\Gamma$ which are $\beta-$preserving diffeomorphisms i.e., induce a commutative diagram

\vspace{2 mm}
\begin{equation*}
\Gamma~\xrightarrow{~\psi~}~\Gamma
\end{equation*}
\begin{equation*}
\beta\downarrow\hspace{12 mm}\downarrow\beta
\end{equation*}
\begin{equation*}
\Gamma_0\xrightarrow{~f~}\Gamma_0
\end{equation*}

\vspace{4 mm}
\noindent
The group of left translations $\psi$ of $\Gamma$ identifies canonically with the group of global sections $\Gamma_a(\Gamma_0,\Gamma),$ the correspondence being given by
\begin{equation}
\psi(g)=\sigma(\epsilon)\cdot g,\hspace{5 mm} \epsilon=\beta(g).
\end{equation}
We can also define \textit{local} right translations whose domains are $\alpha-$saturated open sets of $\Gamma$ and local 
left translations whose domains are $\beta-$saturated open sets. The formula (4) (resp. (5)) establishes an isomorphism 
of the groupoid $\Gamma_{a,loc}(\Gamma_0,\Gamma)$ with the groupoid of local right (resp. left) translations of $\Gamma.$ 

\vspace{3 mm}
\noindent
Given a local 1-parameter family of sections $\sigma_t\in\Gamma(U,\Gamma),$ $|t-t_0|<\epsilon,$ and assuming that $\sigma_{t_0}$ is admissible, then $\phi_*(\frac{d\sigma_t}{dt}|_{t=t_0})\in\Gamma(f(U),\mathcal{L})$ where $f=\beta\circ\sigma_{t_0}$ and $\phi$ is the right translation defined on $\alpha^{-1}(U)$ by $\sigma=\sigma_{t_0}$ (\textit{cf.} (4)). We shall denote this section of $\mathcal{L}$ by $\frac{d\sigma_t}{dt}|_{t=t_0}.$ In particular, the property (1) of the Exponential map implies that
\begin{equation}
\frac{d}{dt}Exp\!~t\Xi|_{t=t_0}=\Xi.
\end{equation}
Let $\Xi\in\Gamma(\Gamma_0,\mathcal{L})$ be global and let $\xi$ be the corresponding right-invariant vector field on $\Gamma.$ Then $Exp\!~t\Xi=(exp\!~t\xi)\circ\iota$ and, since $\xi$ is right-invariant, it follows that $exp\!~t\xi(g)=[Exp\!~t\Xi(\epsilon)]\cdot g$ where $\epsilon=\beta(g).$ We infer that $(exp\!~t\xi)$ is the 1-parameter group of left translations of $\Gamma$ associated to the 1-parameter group of sections $(Exp\!~t\Xi)$ in $\Gamma_a(\Gamma_0,\Gamma).$ Each of these left translations preserves the right-invariant vector field $\xi.$ The uniqueness of the Exponential map is then a consequence of the uniqueness of the 1-parameter group associated to $\xi.$

\vspace{3 mm}
\noindent
When $[\Xi_1,\Xi_2]=0$ then $[\xi_1,\xi_2]=0$ and the 1-parameter groups of left translations commute. It follows that $Exp\!~t\Xi_1\cdot Exp\!~u\Xi_2=Exp\!~u\Xi_2\cdot Exp\!~t\Xi_1=Exp\!~(t\Xi_1+u\Xi_2).$ More generally, we can also prove a Campbell-Hausdorff formula which gives a second order approximation of $Exp\!~t\Xi_1\cdot Exp\!~t\Xi_2$ in terms of $\Xi_1+\Xi_2$ and $[\Xi_1,\Xi_2].$

\vspace{3 mm}
\noindent
Furthermore, we can also consider the 1-parameter group $(\phi_t)$ of right translations associated to $Exp\!~t\Xi$ namely,
\begin{equation}
\phi_t(g)=g\cdot [Exp\!~t\Xi(e)]^{-1}, \hspace{5 mm} e=\alpha(g).
\end{equation}
The vector field $\eta$ associated to $(\phi_t)$ is the left-invariant vector field on $\Gamma$ obtained from $\xi$ by the inversion $g\longmapsto g^{-1}$ of $\Gamma.$ When $\Xi$ is not global, then similar results still hold locally.

\vspace{3 mm}
\noindent
It is not true, in general, that any admissible section $\sigma\in\Gamma_a(\Gamma_0,\Gamma)$ which is \textit{close} to the units section $\iota$ (in any reasonable topology) is of the form $\sigma=Exp\!~t\Xi$ with $\Xi\in\Gamma_c(\Gamma_0,\mathcal{L}).$ The result is not even true locally \textit{i.e.}, in a neighborhood of a point $e\in\Gamma_0.$ In fact, if $\sigma=Exp\!~t\Xi,$ for a given \textit{t}, then the diffeomorphism $\beta\circ\sigma$ is equal to $exp\!~t(\beta_*\circ\Xi).$ There are even local diffeomorphisms arbitrarily close to the identity map and, nevertheless, are not part of the flow of some vector field. It can however be proved, quite easily, that the set $\{Exp\!~t\Xi(e)\!~|\!~\Xi\in\Gamma_c(\Gamma_0,\mathcal{L}),~e\in\Gamma_0,~t\in\textbf{R}\}$ is an open neighborhood of $\Gamma_0$ in $\Gamma.$ Furthermore, taking a local basis of $\mathcal{L}$ in a neighborhood of a point $e\in\Gamma_0,$ we can also define exponential coordinates in a neighborhood of \textit{e} in $\alpha^{-1}(e).$

\vspace{3 mm}
\noindent
We observe that when $\Xi\in\Gamma_c(\Gamma_0,\mathcal{L})\simeq\Gamma_c(\Gamma_0,V\Gamma|_{\Gamma_0})$ takes its values in the Lie algebras of the isotropy groups of $\Gamma$ (i.e., in $T_e\Gamma_e,~e\in\Gamma_0)$ and assuming that each $\Gamma_e$ is a Lie group, then $Exp\!~t\Xi$ takes its values in the isotropy groups of $\Gamma$ and $Exp\!~t\Xi(e)=exp\!~t[\Xi(e)],$ where the right hand side is the usual exponential map $\mathcal{G}\longrightarrow G$ for Lie groups.

\vspace{3 mm}
\noindent
Let $\sigma\in\Gamma_a(\Gamma_0,\Gamma)$ be an admissible section of the Lie groupoid $\Gamma,$ $\phi$ the right translation defined by $\sigma$ (\textit{cf.} (4)) and $\psi$ the left translation also defined by this section (\textit{cf.} (5)). The composite
\begin{equation*}
\psi\circ\phi=\phi\circ\psi:\Gamma\longrightarrow\Gamma,\hspace{3 mm}g\longmapsto\sigma(\epsilon)\cdot g\cdot\sigma(e)^{-1},\hspace{3 mm}e=\alpha(g),\hspace{3 mm}\epsilon=\beta(g),
\end{equation*}
is an automorphism of $\Gamma$ that induces the identity on $\Gamma_0.$ The maps $\sigma\in\Gamma_a(\Gamma_0,\Gamma)\longmapsto\phi_{\sigma}\circ\psi_{\sigma}\in Aut~\Gamma$ and
\begin{equation*}
Ad:\sigma\in\Gamma_a(\Gamma_0,\Gamma)\longmapsto (\phi_0\circ\psi_0)_*\in Aut\mathcal{L}
\end{equation*}
are groupoid morphisms. If $\Xi,\Sigma\in\Gamma(\Gamma_0,\mathcal{L}),$ then the following formula holds,
\begin{equation}
ad~\Xi(\Sigma)=\frac{d}{dt}~Ad(Exp\!~t\Xi)\Sigma=-[\Xi,\Sigma]
\end{equation}
the first equality being the definition.

\vspace{3 mm}
\noindent
We could also define the Lie algebroid sheaf $\mathcal{H}$ by taking the vector bundle $H\Gamma$ composed by the $\beta-$vertical vectors along $\Gamma_0$ and take the bracket of left-invariant vector fields on $\Gamma.$ There would correspond an Exponential map $\Gamma_c(\Gamma_0,\mathcal{H})\longrightarrow\Gamma_b(\Gamma_0,\Gamma)$ where $\Gamma_b(\Gamma_0,\Gamma)$ is the groupoid of all $\alpha-$admissible sections of the fibration $\beta.$ However, the inversion map $g\longrightarrow g^{-1}$ of $\Gamma$ and its differential transport one Exponential map into the other. If $\Xi\in\Gamma_c(\Gamma_0,\mathcal{L})$ and $\Sigma\in\Gamma_c(\Gamma_0,\mathcal{H})$, then $Exp\!~t\Xi$ can be considered as a 1-parameter group of left translations on $\Gamma$ and $Exp\!~t\Sigma$ a 1-parameter group of right translations. It is then obvious that $Exp\!~t\Xi$ and $Exp\!~t\Sigma$ commute.

\section{Prolongations}
Let $\Gamma$ be a Lie groupoid and denote by $J_k\Gamma$ the manifold of all $k-$jets of local sections with respect to the fibration $\alpha:\Gamma\longrightarrow\Gamma_0.$ The fibered manifold $J_k\Gamma$ admits a natural structure of a differentiable category, $j_k\tau(\epsilon)$ being composable with $j_k\sigma(e)$ if and only if $\epsilon=\beta\circ\sigma(e)$ and the resulting composite is equal to $[j_k\tau(\epsilon)]\cdot[j_k\sigma(e)]=j_k s(e),$ where $s(x)=\tau(\beta\circ\sigma(x))\cdot\sigma(x),$ the right hand side being the composition in $\Gamma.$ The units are the $k-$jets of the unit section $\iota:\Gamma_0\longrightarrow\Gamma$ hence can be identified with the elements of $\Gamma_0.$ The groupoid $\Gamma_k$ of all the invertible elements of $J_k\Gamma$ will be called, according to Ehresmann, the $k-$th order prolongation of $\Gamma$. It is an open subset of $J_k\Gamma$ and, with the induced differentiable structure, becomes a Lie groupoid having $\Gamma_0$ as its space of units. We can easily check that $\Gamma_k$ is the set of all jets $j_k\sigma(e)$ such that $j_1(\beta\circ\sigma)(e)$ is invertible or, equivalently, such that $\sigma$ is an admissible local section in a neighborhood of $e$  and it follows that the groupoid structure of $\Gamma_k$ is equal to the extension, to $k-$jets, of the groupoid structure of $\Gamma_{a,loc}(\Gamma_0,\Gamma)$ or, equivalently, the extension to $k-$jets with source in $\Gamma_0$ of the groupoid structure pertaining to the set of all the local right (or left) translations on $\Gamma$. We denote by $\mathcal{L}_k$ its Lie algebroid. When $k=0$ then $\Gamma_0=\Gamma$ and, of course, we shall not bother about using this index since it already has a different meaning. Let $\underline{\Gamma}$ be the sheaf of germs of local admissible sections of $\alpha:\Gamma\longrightarrow\Gamma_0.$ Then $j_k:\underline{\Gamma}\longrightarrow\underline{\Gamma_k}$ is an injective groupoid morphism.

\vspace{3 mm} 
\noindent
We next consider the Lie algebroid $\mathcal{L}$ of $\Gamma$ namely, the sheaf of germs of local sections of the vector bundle $V\Gamma|_{\Gamma_0}$ and denote by $\underline{J_k(V\Gamma|_{\Gamma_0})}\simeq\mathcal{J}_k\mathcal{L}$ the sheaf of germs of local sections of the bundle $J_k(V\Gamma|_{\Gamma_0})$.

\newtheorem{bracket}[PropositionCounter]{Proposition}
\begin{bracket}
There is a unique structure of Lie algebroid on the sheaf $\mathcal{J}_k\mathcal{L}$ satisfying the relation
\begin{equation}
[gj_k\mu,fj_k\eta]=gfj_k[\mu,\eta]+g[\vartheta(\beta_*\circ\mu)f]j_k\eta-f[\vartheta(\beta_*\circ\eta)g]j_k\mu,
\end{equation}
where $\mu,\eta\in\mathcal{L},$\hspace{3 mm}$g,f\in\mathcal{O}_{\Gamma_0},$\hspace{2 mm}  $\beta_*\circ\eta,\beta_*\circ\mu\in\underline{T\Gamma_0}$ and $[\mu,\eta]$ is the bracket in $\mathcal{L}.$ In particular, the map $j_k:\mathcal{L}\longrightarrow\mathcal{J}_k\mathcal{L}$ is an injective Lie algebroid morphism.
\end{bracket}

\noindent
The uniqueness is obvious since any element in $\mathcal{J}_k\mathcal{L}$ is a linear combination of (right) holonomic $k-$jets with coefficients in $\mathcal{O}_{\Gamma_0}$ (left structure). The existence can be proved by transporting the bracket of $\mathcal{L}_k$ \textit{via} the isomorphism of the next theorem. It can also be proved by extending a method, used in a special case, in \cite{Que1967}. This method consists in defining the bracket of two linear combinations by respecting additivity and the formula of the above  Proposition. We can then checks, using the linear operator $\it{\bf{D}}$ of Spencer (\cite{Kumpera1972},\cite{Que1967}), that the definition is independent of the choice of the linear combinations.

\newtheorem{canonical}[TheoremCounter]{Theorem}
\begin{canonical}
$\mathcal{L}_k$ is canonically iomorphic to $\mathcal{J}_k\mathcal{L}.$
\end{canonical}
\noindent
For a proof, we refer the reader to \cite{Kumpera1975}. In the examples to follow, we shall carry out the proof in the special case when $\Gamma_k=\Pi_k\Gamma_0$ is the groupoid of all invertible $k-$jets on a given manifold.

\vspace{3 mm}
\noindent
When $h\leq k,$ the projection $\rho_h:\Gamma_k\longrightarrow\Gamma_h$ is a surjective morphism of Lie groupoids that becomes the identity on the units sub-manifolds (all identical) and induces a surjective Lie algebroid morphism $\rho_{h*}:\mathcal{L}_k\longrightarrow\mathcal{L}_h,$ this last morphism also identifying with the standard projection $\rho_h:\mathcal{J}_k\mathcal{L}\longrightarrow\mathcal{J}_h\mathcal{L}.$ The morphisms $\rho_h$ commute with the Exponential maps \textit{i.e.},
\begin{equation}
\rho_h\circ Exp\!~t\Xi=Exp\!~t\rho_{h*}\Xi,
\end{equation}
where $\Xi\in\Gamma(\Gamma_0,\mathcal{L}_k)$ and $\rho_0=\beta_k:\Gamma_k\longrightarrow\Gamma$ verifies the same properties.

\section{Examples}
The examples that we shall now examine are not simply illustrations of the subject discussed earlier but are, in fact, the basic settings for several applications to be developed in the sequel and relative to the integration of systems of partial differential equations.

\vspace{4 mm}
\noindent
a) \textit{The groupoids} $M^2~and~\Pi_k M.$

\vspace{2 mm}
Let \textit{M} be a manifold and consider the product manifold $N\times M$ where \textit{N} is simply another copy of \textit{M} (i.e., $N=M$). Contrary to the accepted habit, we shall think of \textit{M} as the horizontal component and of \textit{N} as the vertical one in the product. This could be avoided if, instead of the $\alpha-$fibration, we gave preference to the $\beta-$fibration and, instead of the right-invariant vector fields, we opted for the left-invariant fields. The diagonal $\Delta\subset N\times M$ identifies with \textit{M} via the diagonal inclusion $x\longmapsto(x,x).$ We define a groupoid structure on $N\times M$ by means of the following rule. The element $(y,x)$ is composable with $(z,w)$ if and only if $x=z$ in which case the composite is equal to $(y,w).$
For this Lie groupoid structure on $G=N\times M,$ the set of units is the diagonal (that identifies with \textit{M}), the inverse of $(y,x)$ is equal to $(x,y),$ the map $\alpha$ is the second projection $(y,x)\longmapsto x$ and the map $\beta$ is the first projection. This groupoid also satisfies the following \textit{transitivity} condition namely, the fibre product map $\beta\vee\alpha=Id$ is a surjective map of maximal rank. For any $x\in M,$ the fibre $\alpha^{-1}(x)$ is equal to the \textit{vertical} sub-manifold $\{(y,x)|y\in N\}\simeq N$ and similarly for the $\beta-$fibre. The $\alpha-$vertical bundle $V\Gamma$ is equal to the sub-bundle of $T(N\times M)$ composed by all the vertical vectors \textit{i.e.}, those annihilated by $(\pi_M)_*$ hence $V\Gamma|_{\Gamma_0}=V\Gamma|_{\Delta}$ is canonically isomorphic, via $(\pi_N)_*,$ to $TN=TM.$ It follows that $\mathcal{L}\simeq\underline{TM}$ as an $\mathcal{O}_{\Gamma_0}-$linear sheaf. Given two points $x,z\in M,$ there exists a unique right translation $\Gamma_x\longrightarrow\Gamma_z$ of the $\alpha-$fibres namely produced by the element $(x,z)$ and given by $(y,x)\longmapsto(y,z).$ More generally, the right translations of $\Gamma$ are the $\pi_M-$preserving diffeomorphisms of $N\times M$ that commute with $\pi_N$ \textit{i.e.}, are of the form $f\times Id.$ We infer that the right-invariant vector fields $\xi$ on $\Gamma$ are of the form $(\eta,0)$ where $\eta$ is a vector field on $N=M.$ Since $[(\eta,0),(\mu,0)]=([\eta,\mu],0),$ it follows that $\underline{(\pi_N)_*}:\mathcal{L}\longrightarrow\underline{TM}$ is a Lie algebroid isomorphism hence the bracket [\hspace{2 mm},\hspace{2 mm}] on $\mathcal{L}$ is simply the bracket of germs of vector fields on \textit{M}. From the nature of the right-invariant vector fields, we also infer that $Exp\!~t\Xi:x\longmapsto(exp\!~t\eta(x),x),$ where $\eta=(\pi_N)_*\Xi$ and $x\equiv(x,x)\in M.$

\vspace{3 mm}
\noindent
A local section $\sigma:U\longrightarrow\Gamma$ defined on an open set \textit{U} in \textit{M} and with respect to $\alpha$ is a map of the form $\sigma:x\longmapsto(f(x),x),$ where $f=\beta\circ\sigma:U\longrightarrow M$ is differentiable, hence $\sigma=f\times Id.$ The section $\sigma$ is admissible if and only if \textit{f} is a diffeomorphism of \textit{U} onto the open set $f(U).$ If $\tau=g\times Id$ is composable with $\sigma$ then $\tau\circ\sigma=(g\circ f)\times Id.$ We infer that the map
\begin{equation*}
j_k\sigma(x)\in\Gamma_k\longmapsto j_kf(x)=j_k(\beta\circ\sigma)(x)\in\Pi_kM
\end{equation*}
is a Lie groupoid isomorphism that induces a Lie algebroid isomorphism
\begin{equation*}
\mathcal{L}_k\longrightarrow\mathcal{H}_k
\end{equation*}
where $\Gamma_k$ is the $k-$th prolongation of $\Gamma$ and $\mathcal{H}_k$ is the Lie algebroid of $\Pi_kM,$ the bracket being induced by that of the right-invariant vector fields on $\Pi_kM.$ We shall prove that $\mathcal{H}_k$ is canonically isomorphic to $\underline{J_kT}\simeq\mathcal{J}_kT,$ $T=TM,$ the Lie algebroid structure of the latter being defined in the section 10 of \cite{Kumpera1972}. Any local diffeomorphism $\phi$ of \textit{M}, defined on an open set \textit{U}, can be prolongued to a local diffeomorphism $\phi_k$ defined on $\beta^{-1}(U)\subset\Pi_kM$ by setting $\phi_k(X)=[j_k\phi(y)]\cdot X,~y=\beta(X),$ where $\beta$ is the target projection defined on the $k-$jets. If $\theta$ is a local vector field defined on an open set \textit{U}, then the local 1-parameter group $(\phi_t),$ generated by $\theta,$ prolongs to a local 1-parameter group $(\phi_{k,t})$ defined on $\beta^{-1}(U)$ and determines the vector field $\theta_k=\frac{d}{dt}\phi_{k,t}|_{t=0}$ called the $k-$th prolongation of $\theta.$ For any $k-$jet $X\in\beta^{-1}(U),$ the vector $\theta_{k,X}\in T_X\Pi_kM,$ induced by $\theta_k,$ only depends on $j_k\theta(y),~y=\beta(X).$ If $x=\alpha(X),$ then $\theta_{k,X}$ is tangent to the $\alpha-$fibre $\alpha^{-1}(x)$ and the map
\begin{equation}
j_k\theta(y)\in J_kT\longmapsto\theta_{k,X}\in T_X\alpha^{-1}(x)
\end{equation}
is a linear isomorphism. Moreover, if $Y\in\Pi_kM,~\beta(Y)=x$ and $Z=X\cdot Y,$ then $\theta_{k,Z}=(\phi_Y)_*\theta_{k,X},$ where $\phi_Y$ is the right translation produced by \textit{Y}. We infer that $\theta_k$ is a right-invariant local vector field on $\Pi_kM,$ a fact that could also be derived by observing that $(\phi_{k,t})$ is a local 1-parameter group of left translations. Any section $\Xi\in\Gamma(U,J_kT)$ determines, via the isomorphism (11), a right-invariant vector field $\xi$ on the open set $\beta^{-1}(U)$ and, conversely, any right-invariant vector field on $\beta^{-1}(U)$ defines a section of $J_kT$ over \textit{U}. Passing to germs, we obtain the $\mathcal{O}_M-$linear sheaf isomorphism
\begin{equation}
\underline{J_kT}\longrightarrow\mathcal{H}_k
\end{equation}
and we claim that this isomorphism is also a Lie algebroid morphism.

\vspace{3 mm}
\noindent
To prove this, it is sufficient to show (\textit{cf.} \cite{Kumpera1975}, (10.1)) that the Lie algebroid structure of $\mathcal{H}_k,$ transported to $\underline{J_kT},$ satisfies the relation
\begin{equation}
[fj_k\mu,gj_k\eta]=fgj_k[\mu,\eta]+f[\vartheta(\mu)g]j_k\eta-g[\vartheta(\eta)]j_k\mu
\end{equation}
where $f,g\in\mathcal{O}_M,~\mu,\eta\in\underline{TM}$ and $\vartheta(~)$ is the Lie derivative. In fact, using 1-parameter families instead of 1-parameter groups to define the prolongation of vector fields, it is trivial to check that the prolongation process $\theta\mapsto\theta_k,$ for vector fields, is $\bf{R}-$linear and preserves the brackets. Moreover, (12) is $\mathcal{O}_M-$linear and $\vartheta(\xi)(f\circ\beta)=\vartheta(\zeta)f,$ where \textit{f} is any local function on \textit{M}, $\xi$ is the right-invariant vector field on $\Pi_kM$ derermined by a section $\Xi$ of $J_kT$ and $\zeta=\beta_*\xi=\beta\circ\Xi.$ The equality (13) is now evident.

\vspace{3 mm}
\noindent 
In \cite{Kumpera1975}, section 19, the groupoid $\Gamma_kM$ was defined as the quotient of $\widetilde{Aut}~M^2$ modulo a certain equivalence relation $\sim_k$ or, equivalently, as the sheaf of germs of $\beta-$admissible local sections of $\Pi_kM$ with respect to the source fibration $\alpha.$ Furthermore, the group $\Gamma_a(M,\Gamma_kM)$ of admissible global sections was also defined there and it was shown that this group is canonically isomorphic to $\Gamma_a(M,\Pi_kM).$ The uniqueness property of the Exponential maps, defined in section 3 and in \cite{Kumpera1975}, section 19, implies that these Exponential maps are transported one onto the other via the isomorphism
\begin{equation}
\epsilon_k:\tilde{\mathcal{J}}_kT\longrightarrow\mathcal{J}_kT\simeq\underline{J_kT}
\end{equation}
hence,
\begin{equation*}
\Gamma_c(M,\tilde{\mathcal{J}}_kT)\longrightarrow\Gamma_c(M,\underline{J_kT})\hspace{5 mm}and \hspace{5 mm}\Gamma_a(M,\Gamma_kM)\longrightarrow\Gamma_a(M,\Pi_kM)
\end{equation*}
are equivalences.

\vspace{3 mm}
\noindent
The $\ell-$th prolongation of $\Pi_kM$ is not equal to $\Pi_{k+\ell}M.$ In fact, there is a canonical Lie groupoid inclusion $\Pi_{k+\ell}M\longrightarrow\Pi_{\ell}(\Pi_kM)$ whose image is a regularly embedded sub-manifold. Let $\underline{\Pi_1(\Pi_kM)}$ be the sheaf of germs of local sections of $\alpha:\Pi_1(\Pi_kM)\longrightarrow M.$ The non-linear operators $\mathcal{D},$ $\tilde{\mathcal{D}}$ and $\hat{\mathcal{D}},$ defined in \cite{Kumpera1975}, chapters IV and V, are of order 1. They are the composites of $j_1:\Gamma_{k+1}M\longrightarrow\underline{\Pi_1(\Pi_{k+1}M)}$ with (the extension to the germs of sections of) certain bundle maps whose domains are the total space of the locally trivial bundle $\alpha:\Pi_1(\Pi_{k+1}M)\longrightarrow M.$ The explicit description of these bundle maps is left to the reader.

\vspace{3 mm}
\noindent
b) \textit{Fibrations in groups.}

\vspace{2 mm}
\noindent
Let $\pi:\Gamma\longrightarrow M$ be a fibration such that each fibre $\pi^{-1}(x)=\Gamma_x$ is a Lie group. We define a groupoid structure on $\Gamma$ in the following way: The element \textit{h} is composable with \textit{g} if and only if $\pi(h)=\pi(g)$ and, this being the case, the composite in the fibre $\Gamma_{\pi(h)}$ is equal to \textit{hg}, the product in the group. The units of $\Gamma$ are equal to the units of the groups $\Gamma_x,$ hence the set of unit elements identifies with the base manifold \textit{M}. The projections $\alpha$ and $\beta$ are equal to $\pi$ hence the fibres $\alpha^{-1}(x)$ and $\beta^{-1}(x)$ are both equal to the isotropy group  $\Gamma_x$ at the point \textit{x}. This groupoid is called a \textit{fibration in groups} and $\Gamma$ is a differentiable groupoid if and only if the following two conditions are verified:

\vspace{3 mm}
$~(i)$ The units section $\iota:x\longmapsto e_x$ is differentiable and

\vspace{2 mm}
$(ii)$ The fibre product, with respect to $\pi,$ $\Gamma\times_M\Gamma\longrightarrow\Gamma$ is differentiable.

\vspace{3 mm}
\noindent
The differentiability of the inversion map $g\longmapsto g^{-1}$ is a consequence of the condition $(ii)$ and the implicit functions theorem. Since $\pi$ is a fibration, we infer that $\Gamma,$ being a differentiable groupoid, is also a Lie groupoid. Assuming now that this is the case, the vector bundle $V\Gamma|_M,$ (\textit{M} identifying with the units space of $\Gamma$) is a fibration in Lie algebras, the fibre at the point \textit{x} being equal to the Lie algebra $\mathcal{L}(\Gamma_x)$ of the isotropy group. This fibration is not necessarily locally trivial for the algebraic structures but the map
\begin{equation*}
[\hspace{3 mm},\hspace{2 mm}]:V\Gamma|_M\times_M V\Gamma|_M\longrightarrow V\Gamma|_M
\end{equation*}
is differentiable and $\textbf{R}-$bilinear on each fibre. Since the right-invariant vector fields on $\Gamma$ are simply the $\pi-$vertical vector fields $\xi$ whose restriction to each fibre is right-invariant on the group $\Gamma_x,$ we infer that the bracket on the sheaf $\mathcal{L}=\underline{V\Gamma|_M}$ is the extension, to germs of sections, of the fibre-wise bracket in $V\Gamma|_M.$ Since this bracket is $\bf{R}-$bilinear in restriction to each fibre, it follow that the previous sheaf bracket is $\mathcal{O}_M-$bilinear i.e., $[f\mu,g\eta]=fg[\mu,\eta]$ where $f,g\in\mathcal{O}_M$ and $\mu,\eta\in\mathcal{L}$ (\textit{cf.} (1)). Any section $\Xi\in\Gamma(M,\mathcal{L})=\Gamma(M,V\Gamma|_M)$ is global (i.e., $Exp\!~t\Xi$ is global) and $Exp\!~t\Xi(x)=exp\!~t[\Xi(x)]$ where the right hand side is the exponential map $\mathcal{L}(\Gamma_x)\longrightarrow\Gamma_x.$

\vspace{3 mm}
\noindent
Let $\Gamma_k$ be the $k-$th prolongation of the groipoid $\Gamma$ (fibration in groups) and let $\mathcal{L}_k$ be the corresponding sheaf in Lie algebras. Then $\Gamma_k$ is the manifold of all $k-$jets of local sections $\sigma$ of $\Gamma$ with respect to $\pi$ and its groupoid structure is defined as follows: $j_k\tau(y)$ is composable with $j_k\sigma(x)$ if and only if $y=x$ and the composite is then equal to $j_k(\tau\cdot\sigma)(x)$ where $\tau\cdot\sigma$ is the extension, to local sections, of the composition in $\Gamma.$ The units in $\Gamma_k$ are the $k-$jets of the units sections $\iota:x\longmapsto e_x$, the source and the target maps are equal to the extensions $\alpha_k$ and $\beta_k,$ to $k-$jets, of the corresponding projections in $\Gamma$ and the inverse of $j_k\eta(y)$ is equal to $j_k\zeta$ with $\zeta:y\longmapsto(\eta(y))^{-1}.$
We infer that both projections $\alpha_k,\beta_k:\Gamma_k\longrightarrow M$  
are again fibrations in groups.

\vspace{3 mm} 
\noindent
The sheaf $\mathcal{J}_k(\mathcal{L})$ has a unique Lie algebroid structure (in fact, a fibration in algebras) with the bracket satisfying the relation:
\begin{equation}
[fj_k\mu,gj_k\eta]=fgj_k([\mu,\eta]),
\end{equation}
where $f,g\in\mathcal{O}_M$ and $\xi,\eta\in\mathcal{L}$ (\textit{cf.} (9)). By the theorem of section 5, $\mathcal{J}_k\mathcal{L}$ is canonically isomorphic to the algebroid $\mathcal{L}_k$ and we shall now give a direct proof of this result. Since $\Gamma_k$ is a  fibration in groups, the Lie algebroid structure of $\mathcal{L}_k$ is the extension, to germs of sections, of the fibre-wise structure in the Lie algebra  bundle $V\Gamma_k|_M,$ the fibre at each point x being equal to the isotropy algebra of $\Gamma_k$ at that point. We observe next, that $J_k(V\Gamma|_M)$ is a Lie algebra bundle, the bracket being defined by $[j_k\mu(x),j_k\eta(x)]=j_k[\mu,\eta](x),$ where $[\mu,\eta]$ is the extension, to local sections, of the bracket in $V\Gamma|_M.$ The extension, to germs of sections, of the bracket in $J_k(V\Gamma|_M)$ satisfies (15) hence is equal to the canonical bracket on $\mathcal{J}_k\mathcal{L}.$ In view of these remarks, it will suffice to establish a Lie algebra bundle isomorphism
\begin{equation*}
J_k(V\Gamma|_M)\longrightarrow V\Gamma_k|_M.
\end{equation*}
Let $exp_x:\mathcal{L}(\Gamma_x)\longrightarrow\Gamma_x$ be the exponential map and take an element $j_k\zeta(x)\in J_k(V\Gamma|_M).$ For each fixed $t\in\bf{R},~\zeta$ induces the section
\begin{equation*}
\zeta_t:y\in M\longrightarrow exp_y\!~t\zeta(y)\in(\Gamma_k)_x,
\end{equation*}
hence an element $j_k\zeta_t(x)\in(\Gamma_k)_x.$ The map $t\in$\textbf{R}$~\longmapsto j_k\zeta_t(x)\in(\Gamma_k)_x$ is a 1-parameter sub-group hence induces a vector $v=\frac{d}{dt}j_k\zeta_t(x)|_{t=0}\in(V\Gamma|_M)_x$ that satisfies the relation $j_k\zeta_t(x)=exp\!~tv$ where \textit{exp} is the exponential map of $(\Gamma_k)_x.$ We claim that
\begin{equation}
j_k\zeta(x)\in J_k(V\Gamma|_M)\longmapsto v\in V\Gamma_k|_M
\end{equation}
is the desired Lie algebra bundle isomorphism. Since $\frac{\partial}{\partial t}$ and $\frac{\partial^{|\alpha|}}{\partial x^{\alpha}}$ commute,
the map (16) is well defined and injective and a dimensional argument will show that it is also surjective. This map is clearly linear and it is also differentiable since it transforms holonomic sections into differentiable sections of the image bundle. It remains only to prove that (16) preserves the bracket. In fact, $[j_k\zeta(y),j_k\eta(y)]=j_k[\zeta,\eta](y)$ and $[\zeta,\eta](y)=\frac{d}{dt}[\zeta,\eta]_t(y)|_{t=0}=\frac{d}{dt}\Lambda_t(y)|_{t=0},$ where $\Lambda_t(y)=\eta_{\sqrt t}(y)\cdot\zeta_{\sqrt t}(y)\cdot\eta_{-\sqrt t}(y)\cdot\zeta_{-\sqrt t}(y),$ for $t\geq 0,$ and  $\Lambda_t(y)=[\Lambda_{-t}(y)]^{-1},$ for $t\leq 0$ (right-invariant vector fields on a Lie group generate 1-parameter groups of left translations). Again, $\frac{d}{dt}$ commuting with $\frac{\partial^{|\alpha|}}{\partial x^{\alpha}},$ we infer that $\frac{d}{dt}j_k[\zeta,\eta]_t(x)|_{t=0}=\frac{d}{dt}j_k\Lambda_t(x)|_{t=0}=[v,w]\in V\Gamma_k|_M,$ with $w=\frac{d}{dt}j_k\eta_t(x)|_{t=0},$ since $j_k\Lambda_t(x)=[j_k\eta_{\sqrt t}(x)]\cdot[j_k\zeta_{\sqrt t}(x)]\cdot[j_k\eta_{-\sqrt t}(x)]\cdot[j_k\zeta_{-\sqrt t}(x)],$ for $t\geq 0$, $j_k\Lambda_t(x)=[j_k\Lambda_{-t}(x)]^{-1},$ for $t\leq 0,~j_k\zeta_t(x)=exp\!~tv$ and $j_k\eta_t(x)=exp\!~tw.$

\vspace{3 mm}
\noindent
From the last formulas we infer that, for any  $\Xi\in\Gamma(M,J_k(V\Gamma|_M))\simeq\Gamma(M,\mathcal{L}_k),$
\begin{equation*}
Exp\!~t\Xi(x)=j_k\zeta_t(x)=j_kexp\!~t\zeta(x),\hspace{5 mm}\Xi(x)=j_k\zeta(x)
\end{equation*}
and $exp\!~t\zeta(y)=exp_y\!~t[\zeta(y)].$ In particular, $Exp\!~tj_k\zeta=j_kExp\!~t\zeta.$

\vspace{3 mm}
\noindent
A fibration in groups $\pi:\Gamma\longrightarrow M$ is said to be locally trivial with respect to the group structure when, for any $x\in M,$ there exists an open neighborhood \textit{U} of \textit{x}, a Lie group \textbf{H} and a fibration diffeomorphism $\varphi:\pi^{-1}(U)\longrightarrow U\times\bf{H}$ which, restricted to each fibre, is a Lie group isomorphism $\Gamma_x\longrightarrow\{x\}\times\bf{H}\simeq\bf{H}.$ A locally trivial fibration in groups is clearly a Lie groupoid since the trivial fibration $U\times\bf{H}$ is a Lie groupoid and the map $\varphi$ becomes an isomorphism. The bundle $V\Gamma|_M$ then also becomes locally trivial with respect to the algebraic structure since it is locally isomorphic to $U\times\mathcal{H},$ where $\mathcal{H}$ is the Lie algebra of \textbf{H}. The Exponential map and the prolongations can be studied locally in the trivialized patch $\pi^{-1}U$ or, equivalently, in the trivial groupoid $U\times\bf{H}.$

\vspace{3 mm}
\noindent
Let $\Gamma=M\times\bf{H}$ be trivial. Then $V\Gamma|_M=M\times\mathcal{H},$ $Exp\!~t\Xi(x)\simeq exp\!~t[\Xi(x)],$ with $exp:\mathcal{H}\longrightarrow\textbf{H},$ the prolongued groupoid $\Gamma_k$ is canonically isomorphic to the fibration in groups
\begin{equation*}
J_k(M,\textbf{H})=\{j_kg(x)\!~|\!~g:U\longrightarrow\textbf{H},~U~open~in~M\},
\end{equation*}
where the composition is given by $j_kg(x)\cdot j_kh(x)=j_k(g\cdot h)(x)$ and $\mathcal{L}(\Gamma_k)$ is canonically isomorphic to the sheaf in Lie algebras $\underline{J_k(M,\mathcal{H})}.$ Moreover, $J_k(M,\mathcal{H})$ is a bundle in Lie algebras, the bracket being defined by $[j_k\zeta(x),j_k\eta(x)]=j_k[\zeta,\eta](x).$ The isomorphism is given by
\begin{equation}
j_k\zeta(x)\longmapsto\frac{d}{dt}j_k(exp\!~t\zeta)(x)|_{t=0},
\end{equation}
where $exp\!~t\zeta:M\longrightarrow\textbf{H}$ is defined by $(exp\!~t\zeta)(x)=exp\!~t[\zeta(x)].$ We remark that $\Gamma_k\simeq J_k(M,\textbf{H})$ as well as $J_k(M,\mathcal{H})$ are not necessarily trivial and, in fact, the groupoid $J_k(M,\textbf{H})$ is only locally trivial. If \textit{U} is a coordinate neighborhood of \textit{x} then, given any two points $y,z\in U,$ the map $\varphi_{y,z}:J_k(M,\textbf{H})_y\longrightarrow J_k(M,\textbf{H})_z,$ defined by $j_kg(y)\longmapsto j_k(g\circ\phi)(z),$ where $\phi$ is the translation that carries \textit{z} into \textit{y}, is an isomorphism of Lie groups. The map 
\begin{equation*}
\varphi:\alpha_k^{-1}(U)\longrightarrow U\times J_k(M,\textbf{H})_x,\hspace{3 mm}\varphi(X)=(y,\varphi_{y,x}X),~y=\alpha_k X,
\end{equation*}
is a local trivialization with respect to the group structure. Similarly, $J_k(M,\mathcal{H})$ is a locally trivial bundle in Lie algebras. We infer that when the groupoid $\Gamma$ is locally trivial then so is any prolongation $\Gamma_k$ hence also $V\Gamma_k|_M.$

\vspace{3 mm}
\noindent
c) \textit{Trivial Lie groupoids.}

\vspace{2 mm}
\noindent
Let \textit{M} be a manifold, \textbf{H} a Lie group and consider the product manifold $\Gamma=M\times\textbf{H}\times M.$ We define a groupoid structure on $\Gamma$ as follows: The triple $(y,h,x)$ is composable with $(y',h',x')$ if and only if $x=y',$ the composition being then defined by 
$(y,h,x)\cdot(x,h',x')=(y,hh',x').$ The manifold $\Gamma$ becomes clearly a Lie groupoid, the source map 
\begin{equation*}
\alpha:(y,h,x)\in\Gamma\longmapsto x\in M
\end{equation*}
and the target map 
\begin{equation*}
\beta:(y,h,x)\in\Gamma\longmapsto y\in M
\end{equation*}
are fibrations, the units are the elements of the form $(x,e,x)$ where \textit{e} is the unit element in \textbf{H}, the subset of the units identifying consequently with the manifold \textit{M} and the isotropy, at each point, identifies with the group \textbf{H}. Finally, the inverse of $(y,h.x)$ is equal to $(x,h^{-1},y)$. This groupoid is called the \textit{trivial groupoid} with isotropy \textbf{H} and units space \textit{M}. It further satisfies the following transitivity condition: The fibre product
\begin{equation*}
\beta\vee\alpha:\Gamma\longrightarrow M\times M,\hspace{7 mm}\gamma\longmapsto(\beta(\gamma),\alpha(\gamma))
\end{equation*}
is a fibration, it being equal to $\pi_1\times\pi_3.$ We also observe that the groupoid $\Gamma=N\times M$ of the example (\textit{a}) is the special case where the isotropy group \textbf{H} reduces to its unit element. Let us now compute the bracket in $\mathcal{L}.$ We denote by $\mathcal{H}=T_e\textbf{H}$ the Lie algebra of \textbf{H} whose structure is induced by the bracket of right-invariant vector fields on \textbf{H} and observe that, for any $x\in M,$ the $\alpha-$fibre $\Gamma_{\alpha x}=\alpha^{-1}(x)$ identifies canonically with $M\times\bf{H}.$ It follows that $T_x\Gamma_{\alpha x}$ identifies with $T_xM\times\mathcal{H}$ hence $V\Gamma|_M$ identifies with $TM\times\mathcal{H}.$ A section of $V\Gamma|_M,$ defined on an open set \textit{U}, is consequently given by a pair $(\theta,h)$ where $\theta$ is a vector field defined on \textit{U} and $h:U\longrightarrow\mathcal{H}$ is a function. We infer that $\mathcal{L}$ is isomorphic to the sheaf $\underline{TM}\times_M\underline{\mathcal{H}}$ where $\underline{\mathcal{H}}$ is the sheaf of germs of local functions $M\longrightarrow\mathcal{H}.$ Let $\xi$ be a right-invariant vector field on $\Gamma$ and consider its restriction $\eta$ to the fibre $\alpha^{-1}(x)\simeq M\times\bf{H}.$ We can decompose $\eta$ into the sum $\eta=\theta+\tilde{\eta}$ where $\theta=\beta_*\xi=\pi_{1*}\eta$ is a vector field on \textit{M} (lifted to $M\times\bf{H}$) and $\tilde{\eta}$ is a vertical vector field on $M\times\bf{H}$ \textit{i.e.}, $\pi_{1*}\tilde{\eta}=0.$ The right action of the isotropy group $\Gamma_x$ on $\alpha^{-1}(x)$ is given by $(y,h,x)\cdot(x,g,x)=(y,hg,x)$ hence identifies with the action $(y,h)\cdot g=(y,hg)$ of \textbf{H} on $M\times\bf{H}.$ Since $\xi$ is right-invariant, then $\eta$ is also right-invariant by the action of \textbf{H} or, equivalently, the restriction of $\tilde{\eta}$ to each fibre $\{y\}\times\bf{H}$ is a right-invariant vector field on \textbf{H} and consequently the restriction map $\xi\longmapsto\eta$ is an isomorphism of the algebra of right-invariant vector fields on $\Gamma$ onto the algebra of right-invariant vector fields on the $\textbf{H}-$principal bundle $M\times\bf{H}$ since any such invariant field can be uniquely extended to an invariant field on $\Gamma.$ If $\alpha^{-1}(z)\simeq M\times\bf{H}$ is another fibre, then the corresponding restriction map $\xi\longmapsto\mu$ has entirely similar properties and the right translation $\alpha^{-1}(x)\longrightarrow\alpha^{-1}(z)$ produced by the element $(x,g,z)$ is given by $(y,h,x)\cdot(x,g,z)=(y,hg,z)$ hence identifies to $(y,h)\cdot g=(y,hg)$ when both fibres are identified with the corresponding fibres in $M\times\bf{H}.$
The vector field $\xi$ being right-invariant, it follows that $\mu,$ as a vector field on $M\times\bf{H},$ is the right-translated of $\eta$ by the element $g\in\bf{H}$ hence $\mu=\eta.$ From the above discussion, we infer that the algebra of right-invariant vector fields on $\Gamma$ identifies canonically with the corresponding algebra on $M\times\bf{H},$ the identification being obtained by the restriction map on any $\alpha-$fibre. The latter algebra is equal to the set of all the vector fields $\eta=\theta+\tilde{\eta}$ on $M\times\bf{H},$ with $\theta$ a vector field on \textit{M} and $\tilde{\eta}$ a vertical vector field whose restriction to any fibre is right-invariant on \textbf{H}. If $(\theta,h)\in\Gamma(U,V\Gamma|_M)=\Gamma(U,\mathcal{L}),$ then the corresponding right-invariant vector field $\xi$ defined on $\beta^{-1}(U)$ identifies with the right-invariant field $\eta=\theta+\tilde{\eta}$ on $U\times\textbf{H}$ where, for any $y\in U,$ $\tilde{\eta}|_{\{y\}\times\textbf{H}}$ is the right-invariant field on \textbf{H} associated to $h(y).$ Let $(\theta',h')$ be another section which identifies with $\eta'=\theta'+\tilde{\eta}'.$ Then, 
\begin{equation*}
[\eta,\eta']=[\theta,\theta']+[\tilde{\eta},\tilde{\eta}']+[\tilde{\eta},\theta']+[\theta,\tilde{\eta}'],
\end{equation*}
$\pi_{1*}[\eta,\eta']=[\theta,\theta']$ is a vector field on \textit{U} and $[\tilde{\eta},\tilde{\eta}']$ is the right-invariant vertical field on $U\times\textbf{H}$ associated to $[h,h'],$ where the latter bracket is the extension, to functions, of the bracket in $\mathcal{H}.$ Since $\theta=\theta\times 0$ as a vector field on $U\times\bf{H}$ and since $\tilde{\eta}'$ is vertical, we infer that
\begin{equation*}
[\theta,\tilde{\eta}']=\vartheta(\theta)(\phi_{g*}\circ h')(x)=\phi_{g*}[\vartheta(\theta)h'(x)],
\end{equation*}
where $\phi_{g*}:\mathcal{H}\longrightarrow T_g\bf{H}$ is the differential (tangent map) of the right translation on \textbf{H} produced by the element \textit{g}. It follows that the right-invariant vector field 
$[\eta,\eta']$ corresponds to the section
\begin{equation}
([\theta,\theta'],[h,h'])+\vartheta(\theta)h'-\vartheta(\theta')h)
\end{equation}
of $\Gamma(U,\mathcal{L}),$ where $\vartheta(~)$ denotes the Lie derivative and the bracket $[(\theta,h),(\theta',h')]$ 
is equal to (18).

\vspace{3 mm}
\noindent
Let us now examine the Exponential map. The pre-sheaf $\Gamma_a(M,\Gamma)$ identifies with the set of all pairs $(\varphi,g),$ where $\varphi:M\longrightarrow M$ is a diffeomorphism and $g:M\longrightarrow\bf{H}$ a differentiable map, the group structure being 
defined by
\begin{equation*}
(\varphi,g)\cdot(\varphi',g')=(\varphi\circ\varphi',(g\circ\varphi')\cdot g').
\end{equation*}
If $\Xi\in\Gamma(M,\mathcal{L})$ then $\Xi=(\theta,h)$ where $\theta$ is a vector field on \textit{M} and $h:M\longrightarrow\mathcal{H}$ is a differentiable function. The section $\Xi$ is global (\textit{i.e.}, $Exp\!~t\Xi$ is global) if and only if $\theta$ is global (\textit{i.e.}, its 1-parameter group is defined for all values of \textit{t}). Let $\xi$ be the right-invariant vector field on $\Gamma$ that corresponds to $\Xi,$ $\eta=\theta+\tilde{\eta}$ the associated field on $M\times\textbf{H}$ and let us assume that $\theta$ is global. To compute $Exp\!~t\Xi=(exp\!~t\xi)\circ\iota$ on $\Gamma,$ we shall determine the curve $t\longmapsto Exp\!~t\Xi(x)$ for any unit $x\in M$ hence, via the identification $\alpha^{-1}(x)\simeq M\times\textbf{H},$ the integral curve $\gamma(\cdot,x):\textbf{R}\longrightarrow M\times\textbf{H}$ of $\eta$ with initial data $\gamma(0,x)=(x,e).$ Since $\pi_{1*}\tilde{\eta}=0,$ then $exp\!~t\eta$ is $\pi_1-$projectable onto $exp\!~t\theta$ and $\gamma(t,x)=(exp\!~t\theta(x),g(t,x)),$ where $g(\cdot,x):\textbf{R}\longrightarrow\textbf{H}$ satisfies the initial condition $g(0,x)=e$ and
\begin{equation*}
\theta_{exp\!~t\theta(x)}+\tilde{\eta}_{\gamma(t,x)}=\eta_{\gamma(t,x)}=\frac{d\gamma}{dt}(t,x)=\theta_{exp\!~t\theta(x)}+\frac{dg}{dt}(t,x).
\end{equation*}
Since $\tilde{\eta}_{\gamma(t,x)}=\phi_{g(t,x)*}[h\circ exp\!~t\theta(x)],$ (and recalling that $\phi_\lambda$ is the right translation produced by the element $\lambda$) we infer that $g(t,x)$ is the unique solution of the ordinary differential equation
\begin{equation}
\frac{dg}{dt}(t,x)=\phi_{g(t,x)*}\circ h\circ exp\!~t\theta(x)
\end{equation}
with initial data $g(0,x)=e.$

\vspace{3 mm}
\noindent
Given any vector field $\theta$ on a manifold \textit{M}, a Lie group \textbf{H} and a function $h:M\longrightarrow\mathcal{H},$ we can always define the ordinary differential equation (19) which is invariant under the right translations of \textbf{H}. Its behavior is analogous to that of linear equations in the sense that, for any $x\in M,$ the maximal solution $g(t,x)$ with initial data $g(0,x)=g_0$ is defined in the same interval as the right-hand side of the equation namely, the domain of the maximal integral curve of $\theta$ which starts at \textit{x} for $t=0.$ In particular, the equation has global solutions if and only if $\theta$ has global solutions. We now assume that $\textbf{H}\subset GL(n,\textbf{R})$ is a linear Lie group and denote by $\mathcal{H}\subset\mathcal{GL}(n,\textbf{R})$ its linear Lie algebra. The equation (19) reduces in this case to the linear equation
\begin{equation}
\frac{dg}{dt}(t,x)=[h\circ exp\!~t\theta(x)]\cdot g(t,x)
\end{equation}
where, in the right-hand side, the \textit{dot} represents matrix multiplication.

\vspace{3 mm}
\noindent
Let $\Gamma_k$ be the $k-$th prolongation of $\Gamma.$ The elements of $\Gamma_k$ are the $k-$jets $j_k\sigma(x),$ where $\sigma$ is a local admissible section of $\alpha,$ and therefore $\sigma$ identifies with a pair $(\varphi,g)$ where $\varphi$ is a local diffeomorphism of \textit{M} and \textit{g} is a local function with values in \textbf{H} whose domain is equal to the domain of $\varphi.$ The jet $j_k\sigma(x)$ identifies with the pair
\begin{equation*}
(j_k\varphi)(x),j_k g(x))\in\Pi_kM\times_M J_k(M,\textbf{H}),
\end{equation*}
where $J_k(M,\bf{H})$ is the manifold of $k-$jets of local maps $g:U\longrightarrow\bf{H},$ \textit{U} is an open set in \textit{M} and the fibre product is taken with respect to the $\alpha_k$ projections. The composition in $\Gamma_k$ is defined by
\begin{equation}
(j_k\psi(y),j_kh(y))\cdot(j_k\varphi(x),j_kg(x))=(j_k(\psi\circ\varphi)(x),j_k[(h\circ\varphi)\cdot g](x)).
\end{equation}
Recall (example \textit{b}) that $J_k(M,\textbf{H})$ is a Lie groupoid, its Lie algebroid sheaf being isomorphic to $\underline{J_k(M,\mathcal{H})}$ and the groupoid $\Pi_kM$ operates contravariantly (to the right) on $J_k(M,\textbf{H})$ by
\begin{equation*}
(j_kh(y))\cdot(j_k\varphi(x))=j_k(h\circ\varphi)(x),\hspace{5 mm}y=\varphi(x).
\end{equation*}
For any $A\in\Pi_kM,$ the map
\begin{equation*}
\gamma\in J_k(M,\textbf{H})_y\longmapsto\gamma\cdot A\in J_k(M,\textbf{H})_x,\hspace{4 mm}x=\alpha_k(A),\hspace{4 mm}y=\beta_k(A),
\end{equation*}
is an isomorphism of Lie groups and the composition (21) can be re-written by
\begin{equation}
(A',\gamma')\cdot(A,\gamma)=(A'\cdot A,(\gamma'\cdot A)\cdot\gamma),
\end{equation}
where the \textit{dots} indicate the composition of $k-$jets. We infer that
\begin{equation*}
\Gamma_k\simeq\Pi_kM\times_M J_k(M,\textbf{H})
\end{equation*}
is a semi-direct product, where $\Pi_kM$ and $J_k(M,\bf{H})$ are Lie sub-groupoids the latter being normal in the sense that, for any $(Id,\gamma')\equiv\gamma'\in J_k(M,\textbf{H})$ and any $(A,\gamma)\in\Gamma_k$ with $\alpha(A,\gamma)=\beta(\gamma')$ $(=\alpha(\gamma')=\alpha_k(\gamma'),$ the following pertinence $(A,\gamma)\cdot\gamma'\cdot(A,\gamma)^{-1}\in J_k(M,\textbf{H})$ holds. In fact, 
$(A,\gamma)^{-1}=(A^{-1},\gamma^{-1}\cdot A^{-1})$ and consequently
\begin{equation*}
(A,\gamma)\cdot\gamma'\cdot(A,\gamma)^{-1}=(Id,(\gamma\cdot\gamma'\cdot\gamma^{-1})\cdot A^{-1})\equiv(\gamma\cdot\gamma'\cdot\gamma^{-1})\cdot A^{-1}.
\end{equation*}
In particular, if $A=(A,e),$ then
\begin{equation}
A\cdot\gamma'\cdot A^{-1}=\gamma'\cdot A^{-1}.
\end{equation}
Similarly, the Lie algebroid sheaf $\mathcal{L}_k\simeq\mathcal{J}_k\mathcal{L}$ is canonically isomorphic to the fibre product
\begin{equation*}
\mathcal{J}_kT\times_M\underline{J_k(M,\mathcal{H})}\simeq\underline{J_kT}\times_M\underline{J_k(M,\mathcal{H})},
\end{equation*}
the Lie algebroid structure being a semi-direct product of the sub-algebroid $\mathcal{J}_kT,$ with its usual bracket, and the ideal $\underline{J_k(M,\mathcal{H})}$ with the bracket defined at the end of example (\textit{b}). To compute explicitly the bracket in $\mathcal{L}_k,$ we shall first compute the adjoint representation of $\mathcal{J}_kT$ on the ideal $\underline{J_k(M,\mathcal{H})}.$ Observe that $\Pi_kM$ also operates contravariantly on the Lie algebra bundle $J_k(M,\mathcal{H})$ by $A*\chi=\chi\cdot A=j_k(\lambda\circ\varphi)(x)$ where $A=j_k\varphi(x)\in\Pi_kM$ and $\chi=j_k\lambda(y)\in J_k(M,\mathcal{H})$ with $y=\varphi(x).$ For any \textit{A}, the map
\begin{equation*}
\chi\in J_k(M,\mathcal{H})_y\longmapsto A*\chi\in J_k(M,\mathcal{H})_x
\end{equation*}
is the Lie algebra $A-$action induced by the group action $\gamma\longmapsto\gamma\cdot A,$ and we claim that the aforementioned \textit{adjoint} representation is essentially the infinitesimal form of the latter action. Let $\Xi$ be a local section of $J_kT$ and $\Lambda$ a local section of $J_k(M,\mathcal{H}),$ both defined on an open set \textit{U}, and take a local 1-parameter family $(\psi_t)$ of admissible local sections of $\alpha_k:\Pi_kM\longrightarrow M$ such that $\psi_0$ is the identity section and $\frac{d}{dt}\psi_t|_{t=0}=\Xi,$ for example $\psi_t=Exp\!~t\Xi.$ Then the Lie derivative $\vartheta(\Xi)\Lambda=\frac{d}{dt}(\psi_*\circ\Lambda)|_{t=0}$ is a local section of $J_k(M,\mathcal{H})$ defined on \textit{U} and, passing to germs, we obtain a Lie derivative
\begin{equation}
(\Xi,\Lambda)\in \underline{J_kT}\times_M\underline{J_k(M,\mathcal{H})}\longmapsto\vartheta(\Xi)\Lambda\in\underline{J_k(M,\mathcal{H})}
\end{equation}
that is $textbf{R}-$bilinear and satisfies the properties
\begin{equation*}
(a)\hspace{5 mm}\vartheta(\Xi)[\Lambda,\Lambda']=[\vartheta(\Xi)\Lambda,\Lambda']+[\Lambda,\vartheta(\Xi)\Lambda'],
\end{equation*}
\begin{equation}
(b)\hspace{5 mm}\vartheta([\Xi,\Xi'])=[\vartheta(\Xi),\vartheta(\Xi')],\hspace{24 mm}
\end{equation}
\begin{equation*}
(c)\hspace{5 mm}\vartheta(\Xi)(f\cdot\Lambda)=[\vartheta(\Xi)f]\cdot\Lambda+f\cdot\vartheta(\Xi)\Lambda,
\end{equation*}

\vspace{5 mm}
\noindent
where $f\in\underline{J_k}\simeq\mathcal{J}_k$ (germs of $k-$jet of functions), $f\cdot\Lambda$ is the standard $\underline{J_k}-$module structure of $\underline{J_k(M,\mathcal{H})}=\underline{J_kE},$ with $E=M\times\mathcal{H}\longrightarrow M$ and where $\vartheta(\Xi)f$ is defined in the same way as $\vartheta(\Xi)\Lambda$ since $\Pi_kM$ operates on $J_k=J_k(M\times\textbf{R})\simeq J_k(M,\textbf{R}).$ We claim that
\begin{equation}
[\Xi,\Lambda]=\vartheta(\Xi)\Lambda.
\end{equation}
A direct proof using (23), the exponentials of right-invariant vector fields on $\Gamma_k$ and an \textit{adjoint representation} argument is rather strenuous. Quite to the contrary, it is trivial to check that the expression
\begin{equation}
[(\Xi,\Lambda),(\Xi',\Lambda')]=([\Xi,\Xi'],[\Lambda,\Lambda'])+\vartheta(\Xi)\Lambda'-\vartheta(\Xi')\Lambda
\end{equation}
defines a Lie algebroid bracket on $\mathcal{J}_k\mathcal{L}$ that satisfies the formula indicated in the Proposition 1 of section 5, hence it is equal to the canonical bracket of $\mathcal{J}_k\mathcal{L}$ and, in particular, (26) holds.

\vspace{3 mm}
\noindent
\textbf{Remark.} We would expect a minus sign in (26) on account of the inverse in the right-hand side of (23). However, the formula is correct since we are dealing with right-invariant vector fields.

\vspace{3 mm}
\noindent
The relation (27) shows that the bracket in $\mathcal{L}_k$ is formally analogous to (18) even though $\Gamma_k$ is not necessarily a trivial Lie groupoid with base space $(\Gamma_k)_0\equiv M$ and has eventually non-trivial isotropy groups. However, $\Gamma_k$ is locally trivial in a sense that will be specified later. The isotropy group $(\Gamma_k)_x,$ at the point \textit{x}, is isomorphic to $(\Pi_kM)_x\times J_k(M,\textbf{H})_x,$ the group structure being the semi-direct product given by (22). The isotropy algebra is isomorphic to the semi-direct product $(J^0_kT)_x\times J_k(M,\mathcal{H})_x$ given by a formula analogous to (26), where germs of sections of $k-$jets are replaced by the $k-$jets at the point \textit{x} and $(J^0_kT)_x$ is the kernel, at \textit{x}, of the projection $\beta:J_kT\longrightarrow T.$ In fact, $(J^0_kT)_x$ is a Lie algebra for the bracket $[j_k\theta(x),j_k\theta'(x)]=j_k[\theta,\theta'](x),$ the map
\begin{equation*}
j_k\theta(x)\in (J^0_kT)_x\longmapsto\frac{d}{dt}j_kexp\!~t\theta(x)|_{t=0}\in\mathcal{L}(\Pi_kM)_x
\end{equation*}
becomes a Lie algebra isomorphism and, by means of this identification,
\begin{equation*}
exp\!~t[j_k\theta(x)]=j_kexp\!~t\theta(x).
\end{equation*}
The Lie derivative $\vartheta(\Xi)\chi,$ with $\Xi\in(J_k^0T)_x$ and $\chi\in J_k(M,\mathcal{H})_x$ is defined by $\frac{d}{dt}(A_t^*\circ\chi|_{t=0},$ where $t\longmapsto A_t$ is a curve in the isotropy $(\Pi_kM)_x$ such that $A_0=Id$ and $\frac{d}{dt}A_t|_{t=0}=\Xi.$ The proof of the semi-direct decompositions then follow from standard Lie group theoretical arguments but this result can also be proved by simply restricting the section $\Xi$ of $J_kT$ to the bundle of isotropy algebras $J_k^0T$ of $\Pi_kM,$ by calculating in accordance with (24) and by noticing that the value of the section $\vartheta(\Xi)\Lambda,$ at a point \textit{x}, only depends upon the values $\Xi(x)$ and $\Lambda(x)=\chi.$ All the isotropy groups of $\Gamma_k$ are isomorphic with each other and isomorphic to the semi-direct product $(\Pi_k\textbf{R}^n)_0\times J_k(\textbf{R}^n,\textbf{H})_0,$ where $n=dim~M.$ The same result also follows for the isotropy algebras.

\vspace{3 mm}
\noindent
The nature of the bracket (27) enables us to characterize the Exponential map of $\Gamma_k$ or rather the curves $Exp\!~t\Xi(x)$ as the unique solutions, with prescribed initial data, of a system of ordinary differential equations which slightly generalises (19). Nevertheless, we leave the details to the reader since, for practical purposes, it is much easier to trivialize locally the groupoid $\Gamma_k$ and integrate the equation (19).

\vspace{3 mm}
\noindent
\textbf{Remark.} The Lie derivative (24) has a more practical counterpart in the context of \cite{Kumpera1972}, chapter IV. Observing (and also using some of the notations of the above reference) that $J_k(M,\mathcal{H})$ is isomorphic to $J_kE$ with \textit{E} the trivial bundle $M\times\mathcal{H},$ we infer that $\underline{J_k(M,\mathcal{H})}$ is canonically isomorphic to the sheaf  $\underline{M^2\times\mathcal{H}}/\mathcal{I}^{k+1}\underline{M^2\times\mathcal{H}}=\mathcal{J}_kE$ (the  \textit{Malgrange diagonal game}). The Lie algebroid bracket of $J_k(M,\mathcal{H})$ is then the quotient of the bracket on $\underline{M^2\times\mathcal{H}},$ the latter being obtained by extending, to germs of sections, the fibre-wise bracket of the trivial Lie algebra bundle $M^2\times\mathcal{H}.$ A section $\Xi$ of $J_kTM$ identifies with a section of $\underline{J_kTM}\simeq\mathcal{J}_kTM$ hence yields, via the isomorphism $\epsilon_k^{-1},$ a section $\tilde{\Xi}$ of $\tilde{\mathcal{J}}_kTM.$ The local 1-parameter family $\psi_t,$ considered previously, identifies with a local 1-parameter family of admissible sections of $\Gamma_kM$ and, in this context, $\frac{d}{dt}\psi_t|_{t=0}=\tilde{\Xi}$ (\textit{cf.} \cite{Kumpera1972}, sect. 19). We now lift $(\psi_t)$ to a local 1-parameter family of admissible transformations $(\Psi_t)$ defined on $M^2$ and let \textit{L} be a local section of $M^2\times\mathcal{H}$ that represents $\chi.$ Then the Lie derivative $\frac{d}{dt}(\Psi_t^*L|_{t=0},$ $(\Psi_t^*L=L\circ\Psi_t,$ is a local section of $M^2\times\mathcal{H}$ that factors to the local section $\vartheta(\Xi)\chi$ of $\mathcal{J}_kE$ provided by (24). Equivalently, if $\tilde{\xi}$ is a local admissible vector field on $M^2$ which represents $\tilde{\Xi},$ then
\begin{equation*}
\vartheta(\tilde{\xi})L~mod~\mathcal{I}^{k+1}\underline{M^2\times\mathcal{H}}=\vartheta(\tilde{\Xi})\chi.
\end{equation*}
Passing to germs, we infer that the map
\begin{equation}
\tilde{\Xi}\in\tilde{\mathcal{J}}_kTM\longmapsto\vartheta(\tilde{\Xi})\in Der_a\mathcal{J}_kE
\end{equation}
is an injective Lie algebroid representation, where $Der_a\mathcal{J}_kE$ is the Lie algebroid of germs of admissible derivations (of degree zero) with respect to the bracket of $\mathcal{J}_kE$ and where admissibility means the preservation of $\iota_k(\underline{M\times\mathcal{H}})$ as well as of 
\begin{equation*}
J^h_kE=\mathcal{I}^{h+1}\underline{M^2\times\mathcal{H}}/\mathcal{I}^{k+1}\underline{M^2\times\mathcal{H}}
\end{equation*}
in accordance with \cite{Kumpera1972} sect. 27. With the help of these definitions, the proof of the formula (25) is now obvious.

\vspace{3 mm}
\noindent
d) \textit{Locally trivial groupoids.}

\vspace{2 mm}
\noindent
Let $\Gamma$ be a differentiable groupoid, $\Gamma_0=M$ the sub-manifold of units and $\beta\vee\alpha:\Gamma\longrightarrow M\times M$ the map $g\longmapsto(\beta(g),\alpha(g)).$ This map is often called the \textit{anchor} of $\Gamma$ but we shall not do use this terminology and simply refer to it as the fibre product of the two projections. If \textit{U} is an open set in \textit{M}, then the inverse image $\Gamma_U=(\beta\vee\alpha)^{-1}(U\times U)$ is an open differentiable sub-groupoid of $\Gamma.$ We shall say that $\Gamma$ is locally trivial if, for any $x\in M,$ there exists an open neighborhood \textit{U} of \textit{x} such that $\Gamma_U$ is a trivial groupoid in the sense that there exists a trivial groupoid $N\times\textbf{H}\times N$ and an isomorphism of differentiable groupoids $\varphi_U:\Gamma_U\longrightarrow N\times\textbf{H}\times N.$ When $\Gamma$ is locally trivial, then clearly it becomes a Lie groupoid. The local triviality can be expressed in many equivalent ways of which we list e few. The proof of these equivalences is based on the possibility of defining local sections for any fibration (implicit functions theorem).

\newtheorem{liefunctor}[LemmaCounter]{Lemma}
\begin{liefunctor}
Let $\Gamma$ be a differentiable groupoid. Then the following conditions are equivalent:

\vspace{4 mm}
1. $\Gamma$ is locally trivial.

\vspace{3 mm}
2. The fibre product $\beta\vee\alpha$ (or $\alpha\vee\beta$) is a submersion.

\vspace{3 mm}
3. $\alpha$ is a submersion (i.e., $\Gamma$ is a Lie groupoid) and, for any $x\in M,$ the restriction $\beta:\alpha^{-1}(x)\longrightarrow M$ is a submersion. This condition can also be re-phrased interchanging $\alpha$ with $\beta.$

\vspace{3 mm}
4. For any $x\in M,$ the isotropy group $\Gamma_x$ is a Lie group and there exists an open neighborhood U of x together with a differentiable groupoid isomorphism $\varphi:\Gamma_U\longrightarrow U\times\Gamma_x\times U,$ where the image is the trivial Lie groupoid.

\vspace{3 mm}
5. $\beta\vee\alpha:\Gamma\longrightarrow M\times M$ is a locally trivial bundle with a fibre isomorphic to the isotropy $\Gamma_x$ at a point x and group isomorphic to $\Gamma_x\times\Gamma_x$ operating to the left, on the fibre, by $(g',g)\cdot h=g'\cdot h\cdot g^{-1}$ (assuming M connected).

\vspace{3 mm}
6. $\alpha:\Gamma\longrightarrow M$ is a submersion and, for any $x\in M,$ $\beta:\alpha^{-1}(x)\longrightarrow M$ is a principal $\Gamma_x-$bundle (assuming M connected). We can also permute the roles of $\alpha$ and $\beta.$  
\end{liefunctor}

\vspace{3 mm}
\noindent
When $\Gamma$ is locally trivial then it is transitive on any connected component of \textit{M} in the sense that given any two points $x,y$ belonging to the same component, there exists an element $\gamma\in\Gamma$ such that $\alpha(\gamma)=x$ and $\beta(\gamma)=y,$ the same statement being also true when we interchange the roles of $\alpha$ and $\beta.$ We infer that all the isotropy groups and inasmuch all the $\Gamma_x-$principal bundles $\alpha^{-1}(x)$, at points of the same component, are isomorphic as principal bundles. Let $\gamma\in\Gamma$ with $\alpha(\gamma)=y$ and $\beta(\gamma)=x.$ Then the right translation $\phi_{\gamma}:h\in\alpha^{-1}(x)\longmapsto h\cdot\gamma\in\alpha^{-1}(y)$ is a principal bundle isomorphism with respect to the group isomorphism $h\in\Gamma_x\longmapsto\gamma^{-1}\cdot  h\cdot\gamma\in\Gamma_y.$ Moreover, $\alpha:\Gamma\longrightarrow M$ is a locally trivial fibration in principal bundles.

\vspace{3 mm}
\noindent
We now assume that \textit{M} is connected. Let $\xi$ be a right-invariant vector field on $\Gamma$ and let $\eta$ be, as previously, its restriction to a fibre $\alpha^{-1}(x).$ The map $\xi\longmapsto\eta$ is an isomorphism of the Lie algebra of right-invariant vector fields of $\Gamma$ onto the algebra of $\Gamma_x-$invariant vector fields on the principal bundle $\alpha^{-1}(x).$ We infer that the Lie algebroid structure of $\mathcal{L}$ is entirely determined by the bracket of the invariant vector fields on the above principal bundle. If $\Xi\in\Gamma(U,\mathcal{L})$ and $\eta$ is the corresponding right-invariant vector field on $\alpha^{-1}(x)$ defined on the open set $\beta^{-1}(U)$ then, for any $h\in\beta^{-1}(U)$ with $\beta(h)=y,$ the curve $\gamma(t)=[Exp\!~t\Xi(y)]\cdot h$ is the integral curve of $\eta$ with initial value $\gamma(0)=h,$ hence $Exp\!~t\Xi(y)=\gamma(t)\cdot h^{-1}$ is determined by the integration of $\eta.$ If $\Gamma_k$ denotes the $k-$th prolongation of $\Gamma$ and $\alpha_k:\Gamma_k\longrightarrow M$ is the source map of $\Gamma_k,$ then $\alpha_k^{-1}(x)$ is a principal bundle with group $(\Gamma_k)_x.$ Trivializing locally the map $\alpha:\Gamma\longrightarrow M,$ it is easy to see that $(\Gamma_k)_x$ and $\alpha_k^{-1}(x)$ can be defined, up to an isomorphism, by only using $\Gamma_x$ and $\alpha^{-1}(x).$ The construction process is current in the literature under the labels \textit{prolongation of groups} and \textit{prolongation of frame bundles}, these prolongations not being necessarily holonomic (\textit{i.e.}, can be semi or even  non-holonomic).

\vspace{3 mm}
\noindent
Let $\gamma\in\Gamma$ with $\alpha(\gamma)=y$ and $\beta(\gamma)=z.$ Then the left translation $\psi_{\gamma}:h\in\alpha^{-1}(x)_y\longmapsto\gamma\cdot h\in\alpha^{-1}(x)_z$ is a diffeomorphism of the fibres of these bundles which is compatible with the right action of $\Gamma_x.$ Any such diffeomrphism is produced by a left translation hence $\Gamma$ identifies with the groupoid of all the fibre diffeomorphisms of $\alpha^{-1}(x)$ which are compatible with the action of $\Gamma_x.$ Conversely, if \textit{P} is any principal $\textbf{H}-$bundle, then the groupoid $\Gamma$ of all the fibre diffeomorphisms of \textit{P} compatible with the action of \textbf{H} is locally trivial and locally isomorphic to the trivial groupoid $U\times\textbf{H}\times U$ and, thereafter, the spaces $\alpha^{-1}(x)$ and \textit{P} are isomorphic as principal bundles. If $p\in P,$ with $\pi(p)=x,$ then
\begin{equation*}
\varphi:\gamma\in\alpha^{-1}(x)\longmapsto\gamma(p)\in P
\end{equation*}
is such an isomorphism.  

\vspace{3 mm}
\noindent
The previous discussion shows that the theory of locally trivial groupoids is equivalent to the theory of 
principal bundles, the groupoid version perhaps providing a richer understanding.

\vspace{3 mm} 
\noindent
The $k-$th prolongation of a locally trivial groupoid is again locally trivial. In fact, the map $\rho_0=\beta_k:\Gamma_k\longrightarrow\Gamma$ is a surmersion (surjective submersion) and $\beta_k\vee\alpha_k:\Gamma_k\longrightarrow M\times M$ is the composite $(\beta\vee\alpha)\circ\rho_0.$

\vspace{3 mm} 
\noindent
Alternatively, the $k-$th prolongation of a trivialized patch of $\Gamma$ will again be an open sub-groupoid of $\Gamma_k$ and, in fact, will also be a trivialized patch for the prolongation. In the example (c), we saw that the $k-$th prolongation of a trivial groupoid  becomes the semi-direct fibre product $\Pi_kM\times_M J_k(M,\textbf{H})$ which is locally trivial. A local trivialization is provided by the fibre product of a local trivialization of $\Pi_kM$ (in the sense of groupoids) with a local trivialization of $J_k(M,\textbf{H})$ (in the sense of group fibrations) over the same open set \textit{U} and with the induced semi-direct product structure.

\vspace{3 mm} 
\noindent
The groupoids introduced in the examples (a) and (c) are locally trivial. The fibrations in groups considered in (b) are not always locally trivial, in the sense of groupoids, since $\beta\vee\alpha=\alpha\vee\alpha$. Nevertheless, an appropriate notion of local triviality for these objects is apparent.

\vspace{3 mm} 
\noindent
Most problems relative to locally trivial groupoids can be reduced to corresponding problems for trivial groupoids by restriction to trivialized patches. For example, let $\Xi$ be a section of the Lie algebroid $\mathcal{L}$ of a locally trivial groupoid $\Gamma$ and assume that $\Xi$ is global. Taking the restriction $\Xi_U$ of $\Xi$ to a trivialized patch over the open set \textit{U}, the restricted section needs not be global as a section of the Lie algebroid of the corresponding trivial groupoid. Nevertheless, to determine $Exp\!~t\Xi,$ we can calculate a fragment of the curve $Exp\!~t\Xi(x)$ in the trivialized patch namely, the maximal integral curve $\gamma(t),$ with initial data $\gamma(0)=x,$ of the differential equation (19) where $\Xi$ is replaced by $\Xi_U=(\theta_U,h_U).$ Since $Exp\!~t\Xi=(Exp\!~\frac{t}{n}\Xi)^n,$ the calculation of $Exp\!~t\Xi$ reduces to successive calculations, in local trivializations, with the corresponding equations (19). Similarly, the local structure of $\mathcal{L}$ is determined by the corresponding structure of the Lie algebroid in a trivialized patch and the local structure of $\Gamma_k$ is as well determined by the $k-$th prolongation of the groupoid contained in the trivialized patch.These have been described in the example (\textit{c}). 

\vspace{3 mm}
\noindent
The techniques for trivial groupoids can also be applied to the groupoids which satisfy regularity conditions weaker than local triviality. A Lie groupoid is called $\beta-regular$ if, for any $x\in M,$ the restricted map $\beta:\alpha^{-1}(x)\longrightarrow M$ is locally of constant rank (a sub-immersion). By inversion, the same will be true for $\alpha:\beta^{-1}(x)\longrightarrow M$. This condition implies additional structure on the isotropy groups and on the $\Gamma_x-$spaces $\alpha^{-1}(x).$

\newtheorem{nofunctor}[LemmaCounter]{Lemma}
\begin{nofunctor}
Given a Lie groupoid $\Gamma,$ the following conditions are equivalent:

\vspace{4 mm}
1. $\Gamma$ is $\beta-$regular.

\vspace{3 mm}
2. For any $x\in M,$ the isotropy group $\Gamma_x$ has a Lie group structure and the set $M_x=\beta[\alpha^{-1}(x)]$ carries a differentiable sub-manifold structure of M in such a way that $\alpha^{-1}(x),$ with the manifold structure induced by $\Gamma$, becomes a principal $\Gamma_x-$bundle with base space $M_x.$
\end{nofunctor}

\vspace{3 mm}
\noindent
When $\Gamma$ is $\beta-$regular, then the principal $\Gamma_x-$bundle structure on $\alpha^{-1}(x)$ is unique. 
In fact, the underlying manifold structure of $\Gamma_x$ is regularly embedded in $\alpha^{-1}(x)$ hence in $\Gamma$ 
and, further, the topology determines a Lie group structure. As for the manifold structure of $M_x,$ it is also unique since it is a sub-manifold of \textit{M} whose underlying topology is the quotient topology of $\alpha^{-1}(x).$ 
However, the sub-manifold $M_x$ needs not be regularly embedded. If $y\in M_x$ is another point, then $\Gamma_y$ is isomorphic to $\Gamma_x$ via the conjugation $g\in\Gamma_x\longmapsto  h^{-1}gh\in\Gamma_y$ and $\alpha^{-1}(y)$ is isomorphic to $\alpha^{-1}(x)$ via the principal bundle isomorphism $\gamma\in\alpha^{-1}(x)\longmapsto\gamma\cdot h\in\alpha^{-1}(y),$ where \textit{h} is any element in $\alpha^{-1}(y)$ with $\beta(h)=x$ 
(or $h^{-1}\in\alpha^{-1}(x)$). Let $x\in M$ be fixed, write, for simplicity, $M_x=N$ and define 
$\Gamma_N=\alpha^{-1}(N).$ Clearly $\Gamma_N=\beta^{-1}(N)$ hence $\Gamma_N$ is a sub-groupoid of $\Gamma$ with 
unit space \textit{N} and the same isotropy groups as those of $\Gamma$ (along \textit{N}). Since \textit{N} is a sub-manifold of \textit{M}, then $\Gamma_N$ has a corresponding (unique) structure of sub-manifold of $\Gamma$ in 
such a way that $\Gamma_N$ becomes a $\beta-$regular Lie sub-groupoid of $\Gamma$ that induces on each group 
$\Gamma_x$ the original Lie group structure (or topology). $\Gamma_N$ is then, \textit{in abstractio}, a locally
trivial Lie groupoid with units space \textit{N} and isotropy isomorphic to $\Gamma_x.$

\vspace{3 mm}
\noindent
We infer that the isotropy groups of $\Gamma,$ at the points of \textit{N}, form a locally trivial fibration in groups with base space \textit{N} and the restriction $\alpha:\Gamma_N\longrightarrow N$ is a locally trivial fibration in principal bundles. The manifold \textit{M} is the disjoint union of the sub-manifolds $M_x$ whose dimensions need not be locally constant hence they do not form, in general, a foliation of \textit{M}. If \textit{N} and \textit{N'} are two such sub-manifolds, then the global structure of $\Gamma$ induces relations between $\Gamma_N$ and $\Gamma_{N'}$ which nevertheless still remain rather obscure. 

\vspace{3 mm}
\noindent
The $k-$th prolongaion of a $\beta-$regular Lie groupoid is again $\beta-$regular and the restriction map $j_k\sigma(x)\longmapsto j_k(\sigma|N)(x)$ yields a surjective Lie groupoid morphism $(\Gamma_k)_N\longrightarrow(\Gamma_N)_k$ which,in general, is not an isomorphism.

\vspace{3 mm}
\noindent
e) \textit{Linear Lie groupoids.}

\vspace{2 mm}
\noindent
Let $p:E\longrightarrow M$ be a vector bundle and denote by $GL(E)$ the groupoid of all linear isomorphisms between the fibres of \textit{E}. A generic element $\gamma\in GL(E)$ is a linear isomorphism $\gamma:E_x\longrightarrow E_y.$ We set $\alpha(\gamma)=x,~\beta(\gamma)=y$ and the internal operation in $GL(E)$ is simply the composition of maps. There is a unique structure of differentiable groupoid on $GL(E)$ such that, for any trivialization $\phi_U:E_U\longrightarrow U\times\textbf{R}^m,$ the groupoid $GL(E_U)$ is open in $GL(E)$ and the map
\begin{equation*}
\gamma\in GL(E_U)\longmapsto\phi_U\circ \gamma\circ\phi_U^{-1}\in GL(U\times\textbf{R}^m)=U\times GL(m,\textbf{R})\times U
\end{equation*}
is a differentiable groupoid isomorphism, \textit{m} being an appropriate positive integer (the above composition is performed whenever possible). With this structure, $GL(E)$ becomes a locally trivial Lie groupoid called \textit{the general linear groupoid of the vector bundle E}. The isotropy group at a point \textit{x} is equal to $GL(E_x),$ the linear group of $E_x,$ and the isotropy algebra is, of course, $\mathcal{GL}(E_x)$ ($=L(E_x)$) together with the usual bracket of linear transformation $[\lambda,\lambda']=\lambda'\circ\lambda-\lambda\circ\lambda'$ that is also induced by the bracket of right-invariant vector fields The sheafs in Lie algebras $\mathcal{GL}(E)$ and $\mathcal{GL}(E)_k$ associated to $GL(E)$ and to its $k-$th prolongation $GL(E)_k$ have an interesting analytical interpretation  (\cite{Que1968},\cite{Que1969}).

\vspace{3 mm}
\noindent
We shall now prove that $\mathcal{GL}(E)$ is canonically isomorphic to the sheaf $\mathcal{D}(E)$ of germs of linear differential operators $\delta:\underline{E}\longrightarrow\underline{E}$ satisfying the following condition: There exist a vector field $\theta$ on \textit{M} such that, for any local section \textit{s} of \textit{E} and any local function \textit{f}, the following equality holds.
\begin{equation}
\delta(fs)=f\delta(s)+[\vartheta(\theta)f]s
\end{equation}
where $\vartheta(~)$ denotes the Lie derivative. These operators are clearly of order $\leq 1$ and, since any first order linear differential operator $D:\underline{E}\longrightarrow\underline{F}$ satisfies the relation $D(fs)=fD(s)+\sigma(D,df)s,$ where $\sigma(D,df):\underline{E}\longrightarrow\underline{F}$ is the symbol of $D$ evaluated on the co-vector $df,$ we infer that (29) is equivalent to the relation
\begin{equation}
\sigma(\delta,df)s=[\vartheta(\theta)f]s=<\theta,df>s.
\end{equation}
A simple computation will then show that, if the symbol
\begin{equation*}
\sigma(\delta)\in Hom(T^*\otimes E,E)
\end{equation*}
is viewed as a section of $T\otimes E^*\otimes E,$ the condition (30) reads (is equivalent to) $\sigma(\delta)=\theta\otimes id,$ where \textit{Id} is the \textit{identity} section of $E^*\otimes E.$ Any one of the above mentioned equivalent conditions shows that $\theta$ is uniquely determined by $\delta.$ The bracket of $\mathcal{GL}(E)$ identifies with the commutator $[\delta,\delta']=\delta\circ\delta'-\delta'\circ\delta$ of differential operators and the $\mathcal{O}_M-$module structure of $\mathcal{GL}(E)$ identifies with $(f\delta)s=f(\delta s),~f\in\mathcal{O}_M.$ It is straightforward to verify that
\begin{equation*}
[f\delta,g\delta']=fg[\delta,\delta']+f[\vartheta(\theta)g]\delta'-g[\vartheta(\theta')g]\delta.
\end{equation*}

\vspace{3 mm}
\noindent
Let us take $s\in\Gamma(U,E)$ and let $\Phi\in\Gamma_a(V,GL(E))$ be an admissible local section with $\beta\circ\Phi=f$ and $f(V)\subset U.$ Since $\Phi(x)\in Lis(E_x,E_{f(x)}),$ then $\Phi^*\cdot s=\Phi^{-1}\cdot(s\circ f)\in\Gamma(V,E)$ and for any $\Xi\in\Gamma(U,\mathcal{GL}(E)),$ we define the differential operator $\delta:\Gamma(U,E)\longrightarrow\Gamma(U,E)$ by setting
\begin{equation}
\delta(s)=\vartheta(\Xi)s=\frac{d}{dt}((Exp\!~t\Xi)^*\cdot s)|_{t=0}.
\end{equation}
More explicitly, if $\theta=\beta_*\Xi$ then
\begin{equation}
\delta s(X)=\frac{d}{dt}[Exp\!~t\Xi(x)]^{-1}[s\circ exp\!~t\theta(x)]|_{t=0}.
\end{equation}
The operator $\delta$ clearly satisfies the relation (29), with $\theta=\beta_*\Xi,$ since $\exp\!~t\theta=\beta\circ Exp\!~t\Xi$ and, moreover, the mapping
\begin{equation}
\vartheta:\mathcal{GL}(E)\longrightarrow\mathcal{D}(E)
\end{equation}
is the desired $\mathcal{O}_M-$linear Lie algebroid isomorphism. A standard verification will then shows that $\vartheta$ preserves the bracket, though this property as well as the $\mathcal{O}_M-$linearity can also be verified with the help of the local expression (36).
When $\vartheta(\Xi)=\delta=0,$ then $\delta(fs)=[\vartheta(\theta)f]s=0$ implies that $\theta=0.$ We infer that $\Xi$ is a section of the bundle $\bigcup_{x\in M}~\mathcal{GL}(E_x)$ of the isotropy algebras of $GL(E),$ hence
\begin{equation*}
\delta(s)(x)=\frac{d}{dt}(exp\!~t[\Xi(x)])^{-1}[s(x)]|_{t=0}=-[\Xi(x)][s(x)]=0 
\end{equation*}
and consequently $\Xi(x)=0$ \textit{i.e.}, $\Xi=0.$ Let us prove that $\vartheta$ is surjective. Since $Exp\!~(t+u)\Xi=(Exp\!~t\Xi)\cdot(Exp\!~u\Xi)$ (whenever both sides are defined, $\Xi$ being a local section), we infer that
\begin{equation}
\frac{d}{dt}(Exp\!~t\Xi)^*\cdot s=(Exp\!~t\Xi)^*\cdot\delta(s)=\delta\circ(Exp\!~t\Xi)^*\cdot s,
\end{equation}
where $\delta=\vartheta(\Xi)$ and \textit{s} is any local section of \textit{E}. 
Conversely, given a linear differential operator $\delta\in\Gamma(U,\mathcal{D}E)$ with associated vector field $\theta,$ there exists a locally unique differentiable 1-parameter family $(\Psi_t)$ of admissible sections in $\Gamma_a(U,GL(E)),$ with $\Psi_0=Id,$ the identity section, and $\beta\circ\Psi_t=exp\!~t\theta,$ such that
\begin{equation}
\frac{d}{dt}\Psi_t^*\cdot s=\Psi_t^*\cdot\delta(s),
\end{equation}
for any $s\in\Gamma(U,E).$ Moreover, this family is a local 1-parameter group hence $\Psi_t=Exp\!~t\Xi,$ with $\Xi=\frac{d}{dt}\Psi_t|_{t=0},$ and consequently $\vartheta(\Xi)=\delta$ when the derivative in (35) is taken at $t=0.$ This same assertion will be proved later for each trivial patch $U\times GL(m,\textbf{R})\times U$ of $GL(E).$ We observe that $\Gamma_c(M,\mathcal{GL}(E))$ identifies with the Lie algebroid of those linear first order differential operators $\delta:\Gamma(M,E)\longrightarrow\Gamma(M,E)$ that satisfy (29) for compactly supported vector fields $\theta.$ The Exponential map $Exp\!~t\Xi$ of $GL(E),$ $\Xi\in\Gamma_c(M,\mathcal{GL}(E)),$ is then characterized as the unique global solution of (35) with $Exp\!~0$ equal to the identity section and such that $\delta=\vartheta(\Xi).$ Needless to say that $Exp\!~t\Xi$ verifies (34).

\vspace{3 mm}
\noindent
A local section $\Xi\in\Gamma(U,\mathcal{GL}(E))$ takes values in the bundle of isotropy algebras of $GL(E)$ if and only if $\delta=\vartheta(\Xi)$ is $\mathcal{O}_M-$linear \textbf{i.e.}, an operator of order zero. This being the case, $\delta:\Gamma(U,E)\longrightarrow\Gamma(U,E)$ is the extension, to sections, of a vector bundle mapping $\kappa:E_U\longrightarrow E_U$ and $\kappa(v)=-[\Xi(x)]v,$ where $x=p(v).$ If $\Xi'\in\Gamma(U,\mathcal{GL}(E))$ also takes values in the isotropy bundle then
\begin{equation*}
[\Xi,\Xi'](x)=[\Xi(x),\Xi'(x)]=\Xi'(x)\circ\Xi(x)-\Xi(x)\circ\Xi'(x)
\end{equation*}
(the bracket in $\mathcal{GL}(E_x)$ being defined by right invariant vector fields) hence
\begin{equation*}
[\delta,\delta'](x)=[\lambda_x,\lambda_x']=\lambda_x\circ\lambda_x'-\lambda_x'\circ\lambda_x
\end{equation*}
(the bracket in $\mathcal{GL}(E_x)$ being defined by left invariant vector fields).

\vspace{3 mm}
\noindent
We consider next a local trivialization $GL(E_U)\longrightarrow U\times GL(m,\textbf{R})\times U.$ The groupoid $GL(E_U)$ is then represented by the groupoid $GL(U\times\textbf{R}^m)$ where the element $g\equiv(y,g',x):(U\times\textbf{R}^m)_x\longrightarrow(U\times\textbf{R}^m)_y$ is simply the isomorphism $\textbf{R}^m\longrightarrow\textbf{R}^m$ given by the matrix $g'.$ A local section $\Xi\in\Gamma(U,\mathcal{GL}(E))=\Gamma(U,\mathcal{GL}(E_U))$ identifies with a pair $\Xi\equiv(\theta,h)$ where $\theta=\beta_*\Xi$ is a vector field on \textit{U} and $h:U\longrightarrow\mathcal{GL}(m,\textbf{R})$ is a map with 
values in a matrix algebra. The bracket is then given by
\begin{equation*}
[(\theta,h),(\theta',h')]=([\theta,\theta'],-[h,h']+\vartheta(\theta)h'-\vartheta(\theta')h),
\end{equation*}
with $[h,h']=h\circ h'-h'\circ h$ and
\begin{equation}
\delta(s)=\vartheta(\Xi)s=-h\cdot s+\vartheta(\theta)s,
\end{equation}
where $s:U\longrightarrow\textbf{R}^m$ identifies with a section of $U\times\textbf{R}^m.$ In fact, let $\eta=\theta+\title{\eta}$ be the right-invariant vector field on $U\times GL(m,\textbf{R})$ that corresponds to $\Xi=(\theta,h)$ and recall that $\title{\eta}$ is vertical, its restriction to each fibre $\{y\}\times GL(m,\textbf{R})$ being equal to the right-invariant vector field on $GL(m,\textbf{R})$ associated to $h(y).$ Then
\begin{equation*}
[(Exp\!~t\Xi)^*\cdot s](x)=[Exp\!~t\Xi(x)]^{-1}[s\circ exp\!~t\theta(x)]=g(t,x)^{-1}\cdot[s\circ exp\!~t\theta(x)],
\end{equation*}
where $Exp\!~t\Xi(x)=\gamma(t,x)=(exp\!~t\theta(x),g(t,x))$ is the integral curve of  $\eta$ with initial data $\gamma(0,x)=(x,e),$ \textit{e} being the unit matrix, hence
\begin{equation*}
(\delta s)(x)=\frac{d}{dt}~g(t,x)^{-1}|_{t=0}\cdot s(x)+\frac{d}{dt}~s\circ exp\!~t\theta(x)|_{t=0}=
\end{equation*}
\begin{equation*}
=-h(x)\cdot s(x)+[\vartheta(\theta)s](x),\hspace{22 mm}
\end{equation*}
which proves (36). For each $x\in U,$ the equation (34) transcribes into
\begin{equation}
\frac{d}{dt}~(g(t,x)^{-1}\cdot[s\circ exp\!~t\theta(x)])=\hspace{38 mm}
\end{equation}
\begin{equation*}
\hspace{5 mm}=g(t,x)^{-1}\cdot[(-h\cdot s)\circ exp\!~t\theta(x)+(\vartheta(\theta)s)\circ exp\!~t\theta(x)].
\end{equation*}
The left hand side is equal to
\begin{equation*}
[\frac{d}{dt}~(g(t,x)^{-1}]\cdot[s\circ exp\!~t\theta(x)]+g(t,x)^{-1}\cdot ([\vartheta(\theta)s]\circ exp\!~t\theta(x))
\end{equation*}
where after, the equality (37) is equivalent to
\begin{equation}
[\frac{d}{dt}~(g(t,x)]\cdot g^{-1}(t,x)\cdot[s\circ exp\!~t\theta(x)]=[h\circ exp\!~t\theta(x)]\cdot[s\circ exp\!~t\theta(x)].
\end{equation}
The section \textit{s} being arbitrary, we can cancel out the term $s\circ exp\!~t\theta(x)$ in (38) and the resulting equation is equal to (20). Conversely, if $\delta:C^{\infty}(U,\textbf{R}^m)\longrightarrow C^{\infty}(U,\textbf{R}^m)$ is a first order differential operator satisfying the equality (29) then, writing $s=\sum~f_ie_i,$ we can decompose $\delta$ into the sum $\delta_1+\vartheta(\theta),$ where $\delta_1(s)=\sum~f_i\delta(e_i)$ and where $\vartheta(\theta)$ is the Lie derivative of  $\textbf{R}^m-$valued functions.The operator $\delta_1$ is of order zero since it is linear over the functions and consequently is the extension to $C^{\infty}(U,\textbf{R}^m)$ of the differentiable mapping $-h:U\longrightarrow\mathcal{GL}(m,\textbf{R}),$ where $-h(x):\textbf{R}^m\longrightarrow\textbf{R}^m$ is the linear mapping that transforms the basis $\{e_j\}$ into the family $\{\delta(e_j)(x)\}$. It follows that $\delta=\vartheta(\Xi_U),$ with $\Xi_U=(\theta,h)$ and this proves the surjectivity of (33). Furthermore, if $\delta$ is given and if $\theta$ is the associated vector field, then an entirely similar computation will show that, for any $x\in U,$ the equation (35) also reduces to (20), where $\Psi_t(x)=(exp\!~t\theta(x),g(t,x)).$ This proves the assertion relative to the equation (35), bearing in mind the properties of the solution of (20).

\vspace{3 mm}
\noindent
A groupoid $\Gamma$ is called a linear differentiable groupoid if it is a differentiable sub-groupoid of some $GL(E)$ (not necessarily regularly embedded). Restricting \textit{E}, if necessary, we can always assume that the sub-manifolds of units in $\Gamma$ and in $GL(E)$ coincide. The same definition applies for a linear Lie groupoid. The Lie algebroid $\mathcal{L}$ of a linear Lie groupoid $\Gamma$ identifies with a Lie sub-algebroid of $\mathcal{GL}(E),$ hence with a Lie sub-algebroid of $\mathcal{D}(E),$ the identification being $\mathcal{O}_M-$linear. If $\Xi\in\Gamma(U,\mathcal{L})$ and if we agree that  $\vartheta(\Xi)\equiv\delta,$ then $Exp\!~t\Xi\in\Gamma_a(U,\Gamma)$ is given by the solutions of the equation (35) for this given $\delta$ or, equivalently, by the solutions of (20) with $(\theta,h)\equiv\Xi$ (locally). 
If furthermore $\Gamma$ is locally trivial, then it admits a local trivialization of the form $\psi_U:\Gamma_U\longrightarrow U\times \textbf{H}\times U,$ where 
\textbf{H} is a Lie sub-group of $GL(m,\textbf{R}),$ and there exists a local trivialization $\tilde{\psi}_U:E_U\longrightarrow U\times\textbf{R}^m$ of \textit{E} such that the induced trivialization $GL(E_U)\longrightarrow U\times GL(m,\textbf{R})\times U$ restricted to $\Gamma_U$ is equal to $\psi_U.$ The previous computational techniques can be restricted to $U\times \textbf{H}\times U$ since $\Xi\equiv(\theta,h)\in\Gamma(U,\mathcal{L})$ if and only if $h:U\longrightarrow\mathcal{GL}(m,\textbf{R})$ takes values in the Lie algebra of \textbf{H}.

\vspace{3 mm}
\noindent
When $\Psi\in\Gamma_a(U,GL(E)),~f=\beta\circ\Psi$ and $V=f(U),$ then 
\begin{equation*}
\Psi:\Gamma(U,E)\longrightarrow\Gamma(V,E),
\end{equation*}
with
\begin{equation*}
(\Psi\cdot s)(y)=[\Psi\circ f^{-1}(y)]\cdot[s\circ f^{-1}(y)]  
\end{equation*}
or, equivalently, $\Psi\cdot s=(\Psi^{-1})^*\cdot s.$ We infer that the $k-$th prolongation $GL(E)_k$ is a locally trivial Lie sub-groupoid of $GL(J_kE).$ A dimensional argument will then show that $GL(E)_k$ is a proper sub-groupoid. The Lie algebra sheaf (algebroid) $\mathcal{GL}(E)_k$ associated to $GL(E)_k$ is a sub-algebroid of $\mathcal{GL}(J_kE)$ hence also identifies with a sub-algebroid of $\mathcal{D}(J_kE)$ that is moreover an $\mathcal{O}_M-$sub-module. One proves (\textit{cf}. \cite{Que1969}) that this sub-algebroid is the $\mathcal{O}_M-$linear Lie sub-algebroid generated by the $k-$th prolongations of the operators in $\mathcal{D}(E).$ More precisely, let $\Xi$ be a local section of $\mathcal{GL}(E)$ and let $\delta=\vartheta(\Xi)$ be the corresponding section of $\mathcal{D}(E).$ Denoting by $p_k\delta=j_k\circ\delta:\underline{E}\longrightarrow\underline{J_kE}$ the $k-$th prolongation of $\delta,$ this prolongation determines a first order differential operator $j_k\delta:\underline{J_kE}\longrightarrow\underline{J_kE}$ defined by
\begin{equation*}
j_k\delta(\sum~f_ij_ks_i)=\sum~(f_ij_k(\delta(s_i)+[\vartheta(\theta)f_i]j_ks_i)=
\end{equation*}
\begin{equation*}
\hspace{21 mm}=\sum~(f_ip_k\delta s_i+[\vartheta(\theta)f_i]j_ks_i).
\end{equation*}
We can further check that $j_k\delta$ is a section of $\mathcal{D}(J_kE)$ and, clearly, $\sigma(j_k\delta)=\theta\otimes Id.$ Finally, we can also check that $j_k\Xi\in\mathcal{J}_k\mathcal{GL}(E)\simeq\mathcal{GL}(E)_k$ identifies precisely with the element $j_k\delta\in\mathcal{D}(J_kE).$ Since $j_k[\delta,\delta']=[j_k\delta,j_k\delta']$ and since
\begin{equation*}
[f\Lambda,g\Lambda']=fg[\Lambda,\Lambda']+f[\vartheta(\theta)g]\Lambda'-g[\vartheta(\theta')f]\Lambda,\hspace{5 mm}\Lambda,\Lambda'\in\mathcal{D}(J_kE),
\end{equation*}
we infer that the Lie sub-algebroid $\vartheta[\mathcal{GL}(E)_k]\subset\mathcal{D}(J_kE)$ is equal to the $\mathcal{O}_M-$sub-module generated by $j_k\mathcal{D}(E)$ and that the mapping $j_k\Xi\longmapsto j_k\delta$ extends to a Lie algebroid isomorphism $\mathcal{J}_k\mathcal{GL}(E)\longrightarrow\vartheta[\mathcal{GL}(E)_k].$

\vspace{3 mm}
\noindent
Let $\Gamma\subset GL(E)$ be a linear Lie groupoid. Then its $k-$th prolongation $\Gamma_k$ is again a linear Lie groupoid and $\Gamma_k\subset GL(J_kE).$ Denoting by $\mathcal{L}_k$ its Lie algebroid, then $\vartheta(\mathcal{L}_k)$ is equal to the $\mathcal{O}_M-$linear Lie sub-algebroid of $\mathcal{D}(J_kE)$ generted by the $k-$th prolongations of the operators $\delta\in\vartheta(\mathcal{L})\subset\mathcal{D}(E).$

\vspace{3 mm}
\noindent
As an illustration, we consider the following examples: The groupoid $\Pi_1M$ is equal to $GL(TM)$ hence $(\Pi_1M)_k=GL(TM)_k$ is a Lie sub-groupoid of $GL(J_kTM).$ Let $\Pi_{k+1}M\longrightarrow(\Pi_1M)_k$ be the natural inclusion namely, the injective map $j_{k+1}\varphi(x)\longmapsto j_k(j_1\varphi)(x).$ The groupoid $\Pi_{k+1}M$ is then represented as a Lie sub-groupoid of $GL(J_kTM)$ and the action is the standard one namely, $j_{k+1}\varphi(x)\cdot j_k\theta(x)=j_k(\varphi_*\theta)(y),$ with $y=\varphi(x).$ The trivial groupoid $\Gamma=M\times M$ can be considered as a Lie sub-groupoid of $GL(M\times\textbf{R})$ where $(y,x)$ operates as the identity map $(M\times\textbf{R})_x\longrightarrow(M\times\textbf{R})_y.$ In this case, $\Gamma_k=\Pi_kM$ becomes a Lie sub-groupoid of $GL[J_k(M\times\textbf{R})]=GL(J_k)$ ($J_k$ denoting the set of $k-$jets of functions) and the action of this sub-groupoid is the usual one namely, $j_k\varphi(x)\cdot j_kf(x)=j_k(f\circ\varphi^{-1})(y),$ with $y=\varphi(x)$ and \textit{f} any real valued function. If we replace \textbf{R} by any vector space \textit{V}, then we obtain the standard action of $\Pi_kM$ on the $k-jets$ of $V-$valued functions. \nocite{Ackerman1976},\nocite{Almeida1981},\nocite{Brown2010},\nocite{Cartan1893},\nocite{Cartan1932},\nocite{Cartan1937},\nocite{Cartan1938},\nocite{Cartan1952},\nocite{Goursat1891},\nocite{Goursat1898},\nocite{Goursat1905},\nocite{Goursat1922},\nocite{Kumpera1971},\nocite{Kumpera1991},\nocite{Kumpera1998},\nocite{Kumpera2014},\nocite{Kumpera1985},\nocite{Kumpera1982},\nocite{Kumpera1976},\nocite{Matsushima1954},\nocite{Lie1874},\nocite{Lie1876},\nocite{Lie1877},\nocite{Lie1884},\nocite{Lie1885},\nocite{Lie1895},\nocite{Lie1896}

\bibliographystyle{plain}
\bibliography{references}
\end{document}